\documentclass[12pt]{amsart}
\usepackage{amssymb}
\usepackage{graphicx, tikz}
\usetikzlibrary{shapes}
\usepackage{microtype,booktabs}
\usepackage[noadjust]{cite}
\usepackage{xspace}
\usepackage{longtable}
\usepackage{mathtools}
\usepackage{verbatim}
\usepackage[all]{xy}
\SelectTips{cm}{}
\CompileMatrices
\usepackage{xcolor} % for colored links?
\usepackage[pdfusetitle]{hyperref}
\hypersetup{pdfauthor={Paul E. Gunnells, Mark McConnell, Dan Yasaki}} 
\hypersetup{
    colorlinks,
    linkcolor={red!80!black},
    citecolor={blue!50!black},
    urlcolor={blue!80!black}
}

\numberwithin{equation}{section}
\swapnumbers
\newtheorem{theorem}[subsection]{Theorem}

\theoremstyle{definition}

\newtheorem{example} [subsection] {Example}

\newif\ifshowvc
\showvctrue
\showvcfalse  % uncomment this line to turn off inclusion of vc info

\ifshowvc
\input{vc}
\fi

\catcode`\@=11
\def\@secnumfont{\bfseries}
\catcode`\@12

\newif\iflandscapetable
\landscapetabletrue
\DeclarePairedDelimiter{\set}{\{}{\}}
\DeclarePairedDelimiter{\abs}{\lvert}{\rvert}

\newcommand{\magma}{\texttt{Magma}\xspace}
\newcommand{\Q}{\mathbb{Q}}
\newcommand{\QQ}{\Q}
\newcommand{\ZZ}{\mathbb{Z}}
\newcommand{\C}{\mathbb{C}}
\newcommand{\CC}{\C}
\newcommand{\R}{\mathbb{R}}
\newcommand{\Z}{\mathbb{Z}}
\newcommand{\A}{\mathbb{A}}

\newcommand{\n}{\mathfrak{n}}

\newcommand{\cS}{\mathcal{S}}

\newcommand{\OF}{\mathcal{O}_F}

\newcommand{\p}{\mathfrak{p}}

\newcommand{\m}{\mathfrak{m}}

\newcommand{\q}{\mathfrak{q}}
\newcommand{\bG}{\mathbf{G}}

\newcommand{\w}{\omega}

\renewcommand{\curve}[2]{\href{#1}{\texttt{#2}}}
\newcommand{\field}[2]{\href{#1}{#2}}

\newcommand{\paulcomment}[1]{{\color{blue!80!black} (#1)}}

\newcommand{\T}{\mathbf{T}}
\newcommand{\G}{\mathbf{G}}
\newcommand{\heckepoly}{h}
\newcommand{\mat}[1]{{\begin{bmatrix} #1\end{bmatrix}}}

\newcommand{\uu}{\mathbf{u}}

\newcommand{\del}{\partial }
\DeclareMathOperator{\prim}{pr}
\DeclareMathOperator{\new}{new}
\DeclareMathOperator{\HNF}{HNF}

\DeclareMathOperator{\Frob}{Frob}
\DeclareMathOperator{\sgn}{sgn}
\DeclareMathOperator{\Size}{Size}
\DeclareMathOperator{\SL}{SL}%Special linear
\DeclareMathOperator{\GL}{GL}%General linear
\DeclareMathOperator{\U}{U}%unitary group
\DeclareMathOperator{\Gal}{Gal}%Galois  
\DeclareMathOperator{\Sym}{Sym}%Symmetic 
\DeclareMathOperator{\Res}{R}%restriction of scalars
\DeclareMathOperator{\Norm}{N}%norm 
\DeclareMathOperator{\diag}{diag}%diagonal matrix
\def\BorelSerre{\text{BS}}
\def\Eis{\text{Eis}}

\renewcommand{\thefootnote}{\arabic{footnote}}

\begin{document}

\title[Cohomology of $\GL_{3}$ over the Eisenstein integers]{On the cohomology of congruence subgroups of $\GL_{3}$ over
the Eisenstein integers}

\author[Gunnells]{Paul E. Gunnells}
\address{Department of Mathematics and Statistics\\University of
Massachusetts\\Amherst, MA 01003\\
USA}
\email{gunnells@math.umass.edu}

\author[McConnell]{Mark McConnell}
\address{Princeton University\\ Princeton, New Jersey 08540\\ USA}
\email{markwm@princeton.edu}

\author[Yasaki]{Dan Yasaki}
\address{Department of Mathematics and Statistics\\ 
The University of North Carolina at Greensboro\\Greensboro, NC 27412\\
USA}
\email{d\_yasaki@uncg.edu}

\thanks{
PG wishes to thank the National Science Foundation for support of this
research through the NSF grant DMS-1501832.  DY wishes to thank the
National Security Agency for support through the NSA
grant H98230-15-1-0228 and UNCG for support through the Faculty First
Award.} 

\date{\today}

\begin{abstract}
Let $F$ be the imaginary quadratic field of discriminant $-3$ and
$\OF$ its ring of integers.  Let $\Gamma$ be the arithmetic group
$\GL_3 (\OF)$, and for any ideal $\n \subset \OF$ let $\Gamma_{0}
(\n)$ be the congruence subgroup of level $\n$ consisting of matrices
with bottom row $(0,0,*) \bmod \n$.  In this paper we
compute the cohomology spaces $H^{\nu - 1} (\Gamma_{0} (\n); \C)$ as a
Hecke module for various levels $\n$, where $\nu$ is the virtual
cohomological dimension of $\Gamma$.  This represents the first
attempt at such computations for $\GL_{3}$ over an imaginary quadratic
field, and complements work of Grunewald--Helling--Mennicke \cite{ghm}
and Cremona \cite{cremona}, who computed the cohomology of $\GL_{2}$
over imaginary quadratic fields.  In our results we observe a variety
of phenomena, including cohomology classes that apparently correspond
to nonselfdual cuspforms on $\GL_{3}/F$.
\end{abstract}

\maketitle
\ifshowvc
\let\thefootnote\relax
\footnotetext{Base revision~\GITAbrHash, \GITAuthorDate,
\GITAuthorName.}
\fi

\section{Introduction}\label{s:intro}

\subsection{}
Let $F$ be a CM field with ring of integers $\OF$.  Let $\bG$ be the
Weil restriction $\bG =
\Res_{F/\Q} \GL_{n}$, and let $\Gamma \subset \bG (\Q)$ be an
arithmetic subgroup.  Then according to a theorem of Franke
\cite{franke}, the complex group cohomology $H^{*} (\Gamma; \C )$
provides a geometric incarnation of certain automorphic forms; this
should be thought of as analogous to the classical Eichler--Shimura
isomorphism, which computes the (parabolic) cohomology of congruence
subgroups of $\SL_{2} (\ZZ)$ in terms of holomorphic modular forms.
Moreover, thanks to recent work of Harris--Lan--Taylor--Thorne
\cite{hltt} and Scholze \cite{scholze}, we know that these automorphic
forms correspond to certain Galois representations, in the following
sense.  Suppose $\xi \in H^{*} (\Gamma ; \C)$ is a Hecke eigenclass,
let $p$ be a rational prime, and choose an isomorphism $\overline \QQ_{p} \rightarrow \CC$.
Then there is a continuous semisimple Galois representation $\Gal
(\overline F/ F) \rightarrow \GL_{n} (\overline \QQ_{p})$ such that the Hecke eigenvalues
of $\xi$ at almost all primes $\p \in \OF $ can be computed in terms
of the (inverse)
characteristic polynomial of $\rho (\Frob_{\p})$, where $\Frob_\p$ is
the Frobenius at $\p$.  We refer to \S\ref{s:hecke} for more details.

\subsection{}
Hence the complex cohomology of $\Gamma $ provides a concrete method
to compute certain automorphic forms and to explicitly investigate
certain Galois representations.  The goal of this paper is to carry
out this investigation in the case of $\Gamma \subset \GL_{3} (\OF)$,
where $\OF$ is the ring of integers in the imaginary quadratic field
$F$ of discriminant $-3$.  In particular, for any ideal $\n \subset
\OF$ let $\Gamma_{0} (\n) \subset \GL_{3} (\OF)$ be the subgroup of
matrices with bottom row congruent to $(0,0,*) \bmod \n$.  The top
(potentially) nonvanishing cohomological degree of $\Gamma$ is $\nu =
6$, and it is known that the cohomological degrees where the cuspidal
automorphic forms can appear are $3,4,5$, with the contributions to
any degree essentially equivalent to each other.  Then our main
results are the following:
\begin{itemize}
\item We 
compute the cohomology spaces just below the cohomological dimension,
namely $H^{5} (\Gamma_{0} (\n); \C)$,  for 
levels $\n$ of norm $\Norm(\n) \leq 919$.  
\item For most of the levels $\n$ above, we compute the action of the
Hecke operators on the cohomology for a range of primes $\p$.
\item Finally, for most of the eigenclasses we computed, we determine
their motivic nature. 
\end{itemize}
By the last point, we mean that in most cases we were able to identify
the source of Galois representation that is apparently attached to the
eigenclass.  For example, we find that some eigenclasses
correspond to elliptic curves over $F$ via Eisenstein cohomology, and
that some are attached to Hecke characters over $F$.  The
most interesting classes we find are apparently attached to
nonselfdual automorphic forms over $F$.  These are the first
concretely constructed nonselfdual classes over an imaginary quadratic
field, and indeed over any number field other than $\Q$ (classes for
congruence subgroups of $\GL_3/\Q$ were first constructed by
Ash--Grayson--Green \cite{agg} and later by van Geemen--van der
Kallen--Top--Verberkmoes \cite{top}).

\subsection{} To place this work in context, we recall prior work.
Similar computations for $\GL_{2}$ over imaginary quadratic fields
were done by Grunewald--Helling--Mennicke \cite{ghm}, Cremona
\cite{cremona}, Cremona--Whitley \cite{cremona.whitley}, Bygott
\cite{bygott}, and Lingham \cite{lingham}.  In particular these
authors not only computed cohomology, but also compute the Hecke
action and investigated the motivic nature of the cohomology.  The
first cohomology computations for arithmetic groups in $\GL_{3}$ over
an imaginary quadratic field were done by Staffeldt \cite{staffeldt},
who computed the cohomology of $\GL_{3}$ of the Gaussian integers $\ZZ
[\sqrt{-1}]$.  Later computations were done by two of us (PG and DY)
with Dutour Sikiri\'c--Gangl--Hanke--Sch\"urmann \cite{aimpaper} for
$\GL_{3}$ over a variety of imaginary quadratic fields and for
$\GL_{4}$ over the Gaussian and Eisenstein integers.  However, in
these higher $\Q$-rank papers no Hecke operators were computed.
Indeed, it is the Hecke computations, together with the explicit
nonselfdual eigenclasses, that represent the most novel aspects of the
current paper.

\subsection{}
Here is a guide to the paper.  In \S\ref{s:cohcomp} we outline the
techniques used to carry out the cohomology computation and give an
explicit description of the cell complex.  In \S\ref{s:hecke}, we
describe the techniques used to compute the Hecke action on the
cohomology.  Finally, in \S\ref{s:results} we summarize the results of
the computation.

\section{Cohomology computation}\label{s:cohcomp}

\subsection{}
Let $G = \bG (\R) = \GL_{3}(\C)$ be the
group of real points of $\bG$, and let $K = \U (3) \subset G$ be a
maximal compact subgroup.  Let $\Gamma \subset \GL_{3} (\OF)$ be a
finite-index subgroup.  The group $\Gamma$ acts on the symmetric
space $X = G/K$ by left multiplication.  Since  $X$ is
contractible, we have $H^{*} (\Gamma ; \C) \simeq H^{*} (\Gamma
\backslash X; \C)$. 

Thus one can compute the complex group cohomology of $\Gamma$ by computing
cohomology of the locally symmetric space $\Gamma \backslash X$.  This
is the approach we follow, as in the works \cite{agg, top}.  More
precisely, we a generalization of Voronoi's theory of perfect
quadratic forms allows one to construct an $\GL_{3} (\OF)$-equivariant
polyhedral subdivision of $X$.  Such a subdivision was first
constructed by two of us with Dutour
Sikiri\'c--Gangl--Hanke--Sch\"urmann \cite{aimpaper}; we recall the
setup and results below and refer to \cite{aimpaper} for more details.
Using this polyhedral subdivision we form a chain complex that
computes $H^{*} (\Gamma \backslash X; \C)$.

We remark that although the Hecke operators act on cohomology, they do
not act directly on this chain complex.  Computing the Hecke action
requires different techniques, which we describe in \S \ref{s:hecke}.

\subsection{}\label{ss:q-Lambda}
Let $V$ be the $\R$-vector space of $3\times 3$ Hermitian matrices
over $\C$.  Inside $V$ is the cone $C$ of positive-definite Hermitian
matrices.  The cone $C$ is preserved by homotheties $z\mapsto \lambda
z$, $\lambda \in \R_{>0}$, and the quotient $\pi (C) = C/\R_{>0}$ can
be identified with the symmetric space $X$.  Under this
identification, the left action of $\gamma \in \GL_{3} (\OF)$ on $X$
becomes $x\mapsto \gamma x \gamma^{*}$ on $C$, where the star denotes
Hermitian conjugate.

Let $\Lambda = \OF^{3}$, thought of as column vectors, and let
$\Lambda' = \Lambda \smallsetminus \set{0}$.  We have a map
$q \colon \Lambda \rightarrow V$ given by $q (v) = vv^{*}$.  Note that
$q (v) = q (\varepsilon v)$ for any unit $\varepsilon \in
\OF^{\times}$.  As $v$ ranges over $\Lambda'$,
the points $q (v)$ range over various rank one Hermitian forms in the
closure $\bar C \subset V$.  The convex hull $\Pi$ of these images is
called the \emph{Voronoi polyhedron}.

The group $\GL_{3} (\OF)$ acts on $\Lambda$ by left-multiplication and
acts on $\Pi$.  The cones on the faces of $\Pi$ determine a rational
polyhedral fan $\Sigma$ in the closed cone $\bar C$.  Let
$\tilde{\Sigma} \subset \Sigma$ be the subset of cones that meet $C$
nontrivially.  Then the images $\pi (\sigma)$ of the cones $\sigma \in
\tilde{\Sigma}$ form polyhedral cells in $\pi (C) \simeq X$.  We call
this subdivision of $X$ the \emph{Voronoi decomposition}.  

\subsection{} It is clear that the group $\GL_{3} (\OF)$ acts
cellularly on the Voronoi decomposition.  Moreover, one can prove that
there are only finitely many cells modulo $\GL_{3} (\OF)$.  Let
$\Sigma^{*}_{k} (\Gamma )$ denote a set of representatives modulo
$\Gamma $ of the $k$-dimensional Voronoi cells.  One can use the
$\Sigma^{*}_{k}$ to construct a chain complex $(V_{*} (\Gamma), d)$
computing $H^{*} (\Gamma \backslash X; \C)$, although some care must
be taken since the boundary of a cell in $\tilde{\Sigma}$ may contain
faces lying in $\bar C \smallsetminus C$ (such cells are said to lie
\emph{at infinity}), and since one must consider the action of the
stabilizer subgroups.  Details can be found in \cite[\S 3]{aimpaper},
which in turn relies on \cite{gvs}.  Here we content ourselves with
the following statement (cf.~\cite[Theorem 3.7]{aimpaper}):

\begin{theorem}\label{thm:voronoihomology}
Let $\Gamma \subset \GL_{3} (\OF)$ be a finite index subgroup, and let
$V_{k} (\Gamma )$ be free abelian group on the Voronoi cells
$\Sigma^{*}_{k} (\Gamma )$ mod $\Gamma $.  Then there is a differential $d\colon
\Sigma^{*}_{k}(\Gamma ) \rightarrow \Sigma^{*}_{k-1} (\Gamma )$ such that the homology of
the resulting complex $(V_{*}\otimes \CC, d)$ is isomorphic to the
cohomology $H^{*} (\Gamma \backslash X; \CC )$.  More precisely, we
have an isomorphism
\[
H_{k} (V_{*}\otimes \CC) \longrightarrow H^{8-k} (\Gamma \backslash X; \CC).
\]  
\end{theorem}

\subsection{}\label{ss:vorcones} We finish this section by giving a
complete description 
of the sets $\Sigma^{*}_{k}$ for the case of interest, $F =
\Q(\sqrt{-3})$.  We encode cells by giving the lists of 
vectors in $\Lambda$ that give the vertices of the corresponding faces
of $\Pi$.  Such vectors are called \emph{minimal vectors}, because
they correspond to the vectors on which certain positive-definite
Hermitian forms attain their minima.

Let $\omega
= \frac{1 + \sqrt{-3}}{2} \in \OF$.  Let $A$ denote the matrix
\setcounter{MaxMatrixCols}{20}
\[
A = \begin{bmatrix}
1 & 1 & 1 & 1 & 1 & 1 & 0 & 0 & 0 & 1 & 1 & 1 & 0\\ 
-\omega - 1 & -\omega - 1 & -\omega & -1 & -\omega & -1 & 1 & 1 & 0 & -\omega & -1 & 0 & 1\\
1 & \omega & 0 & 0 & \omega & 1 & -1 & 0 & 1 & \omega^2 & \omega & 0 & -\omega
\end{bmatrix}.
\]
We describe $\Sigma^*$ in terms of the columns of $A$.  In what
follows, the list of integers $[m_1, \dots, m_k]$ denotes the cell in
$X$ that is the image of the cone $\sigma \in \tilde{\Sigma}$
generated by the points $q (v)$, where $v$ is taken from the columns
of $A$ indexed by $m_1, \dots, m_k$.

The top dimensional cells in $\Sigma^{*}$ have dimension that of $\dim
X = 8$.  There are two equivalence classes $8$-dimensional cells, with representatives
\[
a_1 = [ 1, 2, 3, 4, 5, 6, 7, 8, 9 ]\\
\quad \text{and} \quad a_2 = [ 2, 3, 4, 5, 6, 7, 8, 9, 10, 11, 12, 13 ].  
\]
Recall that the $f$-vector of a $d$-dimensional polytope $\sigma$ is
the vector $(f_{-1}, f_0, f_1, \dots, f_{d})$, where $f_k$ is the number of
$k$-dimensional faces of $\sigma$ and $f_{-1} := 1$.  The $f$-vector of $a_1$ is
\[f(a_1) = ( 1, 9, 36, 84, 126, 126, 84, 36, 9, 1 ).\]
Since $a_1$ has nine vertices in an $8$-dimensional space, we see that
$a_1$ is a simplex, and this is indeed compatible with its $f$-vector.
As such, the faces of $a_1$ are easily described as all of the subsets
of $\set{1, 2, \dots, 9}$.  The cone $a_2$ is not a simplex, but it is
a simplicial polytope since all of its $81$ facets are simplices.  The
$f$-vector of $a_2$ is
\[f(a_2) = ( 1, 12, 66, 216, 459, 648, 594, 324, 81, 1 ).\]

The $7$-dimensional cells in $\Sigma^*$ fall into two equivalence
classes under the action of $\Gamma$, namely
\[
  b_1 = [ 1, 2, 3, 4, 5, 6, 7, 8 ] \quad \text{and} \quad 
b_2 =[ 4, 6, 7, 8, 9, 11, 12, 13 ];
\]
these are both simplexes (this implies all lower-dimensional cells
are also simplexes).  All nine of the facets of $a_1$ are of type
$b_1$.  Of the $81$ facets of $a_2$, nine are of type $b_1$ and $72$
are of type $b_2$.

The $6$-dimensional cells in $\Sigma^*$ fall into three equivalence
classes
\[
  c_1 =[ 1, 2, 3, 4, 5, 6, 7 ],\quad 
c_2 = [ 4, 6, 8, 9, 11, 12, 13 ],\quad \text{and} \quad 
c_3 = [ 4, 6, 7, 8, 9, 11, 13 ].
\]
The eight facets of $b_1$ are type $c_1$.  There are six facets of $b_2$
of type $c_2$ and one each of types $c_1$ and $c_3$.

The $5$-dimensional cells of $\Sigma^*$ fall into four equivalence classes
\begin{align*}
d_1 &= [ 1, 2, 3, 4, 5, 6 ],&
d_2 &= [ 1, 2, 3, 4, 5, 7 ],\\
d_3 &= [ 4, 6, 8, 9, 11, 12 ],&
d_4 &= [ 4, 8, 9, 11, 12, 13 ].
\end{align*}
The cell $c_1$ has one facet of type $d_1$ and six facets of type
$d_2$.  The cell $c_2$ has three facets of type $d_2$, three facets of
type $d_3$, and one facet of type $d_4$.  The cell $c_3$ has one facet
of type $d_1$ and three facets of type $d_3$.

The $4$-dimensional cells of $\Sigma^*$ fall into three equivalence classes
\[
e_1 = [ 1, 2, 3, 4, 5 ], \quad
e_2 = [ 1, 2, 3, 4, 7 ],\quad \text{and} \quad 
e_3 = [ 4, 6, 8, 9, 11 ].
\]
All six of the facets of $d_1$ are type $e_1$.  The cell $d_2$ has three
facets of type $e_1$ and $e_2$ each.  The cell $d_3$ has one facet of
type $e_1$, three facets of type $e_2$, and two facets of type $e_3$.
All six of the facets of $d_4$ are type $e_2$.

The $3$-dimensional cells of $\Sigma^*$ fall into two equivalence classes
\[
f_1 = [ 1, 2, 3, 4 ]\quad \text{and} \quad 
f_2 = [ 1, 2, 4, 5 ].
\]
The cell $e_1$ has three facets of type $f_1$ and two facets of type
$f_2$.  The cell $e_2$ has one facet of type $e_1$ and four facets of
type $f_2$.  The cell $e_3$ has one facet at infinity and four facets
of type $f_2$.

The $2$-dimensional cells of $\Sigma^*$ fall into a single equivalence class
\[
g_1 = [1, 2, 3].
\]  
All of the facets of $f_1$ are type $g_1$.  The cell $f_2$ has one
facet at infinity and three facets of type $g_1$.  The boundary of
$g_{1}$ consists of three $1$-cells, all of which lie at infinity.

Since our goal is to compute $H^{5} (\Gamma ; \CC)$, we are primarily
interested in the Voronoi complex in degrees $2, 3, 4$.

\section{Hecke operators}\label{s:hecke}
\subsection{}
In this section we describe the techniques used to compute the Hecke
action on the cohomology.  Our basic technique is that of
\cite{experimental}, adapted to the setting of congruence subgroups
$\Gamma \subset \GL_{3} (\OF)$.  This
technique has already been applied in several other settings:
$\GL_{2}$ over a totally complex quartic field \cite{qzeta5, jones},
$\GL_{2}$ over a real quadratic field \cite{heckants}, 
$\GL_{2}$ over a non real cubic field \cite{neg23paper,complexcubic},
and $\GL_{4}$ over $\Q$ \cite{computation,AGM2,AGM3,AGM4}.  The
common feature all these cases have is that the cohomology group of
interest is $H^{\nu -1}(\Gamma; \C)$, where $\nu$ is the virtual cohomological
dimension of the relevant arithmetic group.  Indeed, for
$\GL_{n}$ over a number field, the cohomological degrees where the
cuspidal automorphic forms can contribute to cohomology (the
\emph{cuspidal range}) depends only on $n$ and the signature of the field.
The above cases, together with the current paper, cover all
cases where the top of the cuspidal range coincides with $\nu -1$.

Before presenting details of our method, we remark that the current
case ($\GL_{3}$ over an imaginary quadratic field $F$) is much more
like $\GL_{4}$ over $\Q $ than the $\GL_{2}$ examples cited above.
This is because the field $F$ has no units of infinite order, which
implies that certain complicated constructions in
\cite{qzeta5,complexcubic,heckants} are not necessary.

\subsection{} We begin by describing the construction of the Hecke
operators.  Fix a level $\n$, let $\Gamma = \Gamma_{0} (\n)$, and let
$\p$ be a prime ideal not dividing $\n$.  Let $\pi$ be a generator for
$\p$ ($\p$ is always principal since $F$ has class number $1$).
Consider the two diagonal matrices $g_{1} = \diag(1, 1, \pi)$ and
$g_{2} = \diag(1, \pi, \pi)$.  The matrices $g_{i}$ each determine
correspondences on $\Gamma \backslash X$ as follows.  Let $\Gamma ' =
g^{-1}\Gamma g \cap \Gamma$, where $g$ either equals $g_{1}$ or
$g_{2}$.  The group $\Gamma '$
has finite index in both $\Gamma$ and $g^{-1}\Gamma g$.  We have a
diagram $C (g)$
\[
\xymatrix{&{\Gamma'\backslash X}\ar[dl]_{s}\ar[dr]^{t}&\\
 {\Gamma\backslash X}&&{\Gamma\backslash X}}
\]
called a \emph{Hecke correspondence}.  The map $s$ is induced by the
inclusion $\Gamma '\subset \Gamma$, while $t$ is induced by the
inclusion $\Gamma '\subset g^{-1}\Gamma g$ followed by the
diffeomorphism $g^{-1}\Gamma g\backslash X\rightarrow \Gamma
\backslash X$ given by left multiplication by $g$.  Specifically, 
\[
s (\Gamma 'x) = \Gamma x, \quad t (\Gamma 'x) = \Gamma gx, \quad x\in X.
\]
The maps $s$ and $t$ are finite-to-one, since the indices $[\Gamma
':\Gamma]$ and $[\Gamma ':g^{-1}\Gamma g]$ are finite.   This implies
that we obtain maps on cohomology 
\[
s^{*}\colon H^{*} (\Gamma \backslash X) \longrightarrow
H^{*} (\Gamma' \backslash X), \quad t_{*}\colon H^{*} (\Gamma' \backslash X) \longrightarrow
H^{*} (\Gamma \backslash X).
\]
Here the map $s^{*}$ is the usual induced map on cohomology, while the
``wrong-way'' map $t_{*}$ is given by summing a class over the finite
fibers of $t$.  These maps can be composed to give a map
\[
T_{g} := t_{*}s^{*}\colon H^{*} (\Gamma \backslash X;
\C )\longrightarrow  H^{*} (\Gamma \backslash X;
\C ),
\]
which is called the \emph{Hecke operator} associated to $g$.  There is
an obvious notion of isomorphism of Hecke correspondences.  One can
show that up to isomorphism, the correspondence $C (g)$ and thus the
Hecke operator $T_{g}$ depend only on the double coset $\Gamma
g\Gamma$.  

\subsection{}\label{ss:cosetreps}
Thus for each prime $\p$ prime to $\n$, we have defined two operators
on cohomology which we denote $T (\p,k)$, $k=1,2$.  To compute them, we need to
understand how the double cosets $\Gamma g_{i} \Gamma$ break up into
left cosets.  More precisely, for $g \in \{g_{1},g_{2} \}$, the double
coset $\Gamma g \Gamma$ can be written as a disjoint union of left
cosets,
\[\Gamma g \Gamma = \coprod_{h \in \Omega} \Gamma h\]
for a certain finite subset $\Omega$ of $3 \times 3$ matrices with
entries in $\OF$.  We now compute $\Omega$ explicitly for the two
operators.  Let $R \subset \OF$ denote a set of representatives for
the residue field $\OF/\p$.  For the operator $T(\p,1)$, the set
$\Omega$ can be taken to be
\[\Omega = \set*{\mat{\pi \\ & 1 \\ && 1}} \cup \set*{ \mat{1 &&
    \alpha\\ & 1 & \beta\\ && \pi}  \colon \alpha, \beta \in R} \cup
\set*{\mat{1 & \alpha &\\ & \pi & 
  \\ && 1} \colon \alpha \in R }.\]
For $T(\p,2)$, the set $\Omega$ can be taken to be
\[\Omega = \set*{\mat{\pi \\ & \pi \\ && 1}} \cup \set*{ \mat{1 &
    \alpha& \beta
    \\ & \pi & \\ && \pi}  \colon \alpha, \beta \in R} \cup
\set*{\mat{\pi &  &\\ & 1 & \alpha
  \\ && \pi} \colon \alpha \in R }.\]

If $f$ is an eigenclass with $T(\p,k) = a(\p,k) f$, then we define the
$\GL_3$ \emph{Hecke polynomial}
\begin{equation}\label{eq:heckepoly}
\heckepoly_{\p} (t) = \heckepoly(f,\p) = 1 - a(\p,1) t + a(\p,2)\Norm(\p)t^2 - \Norm(\p)^{3}t^3.  
\end{equation}  
If we set $t = \Norm(\p)^{-s}$, where $s$ is a complex variable, then
the resulting Dirichlet polynomial is the inverse of the local factor
at $\p$ of the $L$-function attached to $f$.  It is also this
polynomial that, for almost all primes $\p$, matches the inverse
characteristic polynomial of $\rho (\Frob_{\p})$, where $\rho$ is the
Galois representation attached to the eigenclass.

\subsection{} As mentioned before, the Hecke operators do not act on the
Voronoi cells, and thus do not directly act on the Voronoi complex.
To remedy this, we define another complex with Hecke action, the
\emph{sharbly complex}, that also computes the cohomology of $\Gamma$.
This complex was originally introduced by Ash \cite{ash.sharb}.

Recall that for any $x\in \Lambda'$, we have
constructed a point $q (x)\in \bar{C}$ (see
\S\ref{ss:q-Lambda}). 
Write $x\sim y$ if $q (x)$ and $q (y)$ determine the same point in
$\bar{X}$; this happens if and only if the rays generated by $q (x)$
and $q (y)$ coincide, i.e., $\R_{\geq 0} q(x) = \R_{\geq 0} q(y)$.  Let
$\cS_k$, $k\geq 0$, be 
the $\Gamma$-module $A_{k}/C_{k}$, where $A_{k}$ is the set of formal
$\CC$-linear sums of symbols $\uu =[x_1, \dots, x_{k+3}]$, where each
$x_i$ is in $\Lambda'$, and $C_{k}$ is the
submodule generated by
\begin{enumerate}
\item $[x_{\sigma(1)}, \dots, x_{\sigma(k+3)}]-\sgn(\sigma)[x_1,
\dots, x_{k+3}]$,
\item $[x, x_2, \dotsc, x_{k+3}] - [y, x_2, \dotsc , x_{k+3}]$ if $x\sim
y$, and 
\item $\uu $ if $x_1, \dotsc  ,
x_{k+3}$ are contained in a hyperplane (we say $\uu $ is \emph{degenerate}).
\end{enumerate}  
We define a boundary map $\partial\colon \cS_{k+1} \to \cS_{k}$ by
\begin{equation}\label{eq:boundary}
\del [x_1, \dotsc  , x_{k+3}] =\sum_{i=1}^{k+3}(-1)^i[x_1, \dotsc , \hat{x}_i,\dotsc  , x_{k+3}],
\end{equation}
where $\hat{x}_{i}$ means omit $x_{i}$.
The resulting complex $(\cS_{*}, \partial )$ is called the
\emph{sharbly complex}.  

\subsection{}
The group $\Gamma$ acts on $\cS_*$ via the obvious left action of
$\Gamma$ on $\Lambda$, and one can show that the homology of the
complex of coinvariants $(\cS_{*})_{\Gamma}$ is isomorphic to the
cohomology $H^{*} (\Gamma ; \C)$.  More precisely, we have an
isomorphism
\[
H_{i} ((\cS_{*})_{\Gamma}) \rightarrow H^{\nu -i} (\Gamma ; \C),
\]
where $\nu = 6$ is the virtual cohomological dimension of $\Gamma$;
for details see \cite[\S1.4]{ash.sharb}.  The Hecke operators also act on
cohomology through this action.  More precisely, let $\xi \in \cS_*$
be a sharbly chain that is a cycle in $(\cS_{*})_{\Gamma}$, and write
$\xi = \sum n (\uu) \uu $.  Then if $\Omega$ is the set of coset
representatives described in \S\ref{ss:cosetreps} for the Hecke
operator $T$, then the chain
\begin{equation}\label{eq:heckim}
T (\xi) = \sum n (\uu) \sum_{g\in \Omega} g\cdot \uu
\end{equation}
is a well-defined cycle mod $\Gamma$.

\subsection{} Recall that all Voronoi cells of dimension less than or
equal to $4$ are simplices (\S\ref{ss:vorcones}).  Thus in these
dimensions an oriented Voronoi $k$-cell corresponds to an ordered list
of $k+1$ vectors in $\Lambda$.  This means one can naturally identify
subspaces of $\cS_{(k+1)-3}$ for $k=2,3,4$ corresponding to the
Voronoi cells of these dimensions.  This identification is compatible
with the boundary maps and the $\Gamma$-action, so we can think of
$V_{*} (\Gamma )$ as being a subcomplex of $(\cS_{*})_{\Gamma }$ in
these degrees.  The main task in computing the Hecke action on
$H^{5}(\Gamma; \CC)$, our cohomology group of interest, is then the
following: given any cycle $\xi \in (\cS_{1})_{\Gamma}$ representing a
cohomology class in $H^{5} (\Gamma ; \CC)$, find a cycle $\xi ' \in
(\cS_{1})_{\Gamma}$ supported on $V_{3} (\Gamma)$ equivalent to it.

\subsection{} We now describe how to accomplish this.  Given a point
$q (x)$ with $x\in \Lambda$, let $R(x)$ be the point $y\in \Lambda'$
such that $y\sim x$ and $q(y)$ is closest to the origin along the ray
$\R_{\geq 0} q (x)$.  We call $R (x)$ the \emph{spanning point} of
$x$; for any sharbly $\uu = [x_{1},\dots ,x_{n}]$ we can speak of its
set of spanning points.  Any sharbly $\uu = [x_{1},\dots ,x_{n}]$
determines a closed cone $\sigma (\uu)$ in $\bar C$.  We call a
sharbly $\uu = [x_{1},\dots ,x_{n}]$ \emph{reduced} if its spanning
points are a subset of the spanning points of some fixed Voronoi cone;
similarly we call a sharbly cycle reduced if each sharbly in its
support is reduced.  Note that since Voronoi cells are not simplicial
in general, if $\uu$ is reduced it does not mean that $\sigma (\uu)$
is the face of some Voronoi cell, and thus a reduced sharbly cycle
need not come from a Voronoi cycle.  However, it is clear that there
are only finitely many reduced sharbly cycles modulo $\Gamma$.
Moreover, it is not difficult to write any reduced sharbly cycle in
terms of Voronoi cycles directly, so our main challenge is to rewrite
$T (\xi)$ in terms of reduced sharbly cycles.

Computing directly if a sharbly is reduced is very expensive.  Hence
we introduce another measure on cycles, called \emph{size}.  For any
$0$-sharbly $\uu = [x_{1}, x_{2}, x_{3}]$, define the size of $\uu $
to be $\Size (\uu) = \abs{\Norm( \det (R (x_{1}), R (x_{2}),
R(x_{3}))}$.  We extend this to general sharblies by first defining
the size of $\uu = [x_{1},\dotsc ,x_{n}]$ to be the maximum of the
size of the sub $0$-sharblies $[x_{i},x_{j}, x_{k}]$, and then
defining the size $\Size (\xi)$ of a chain $\xi = \sum n(\uu)\uu $ to
be the maximum size occuring in its support.  Clearly the smallest
size any chain can have is $1$.  We remark that any chain of size $1$
is reduced, but the converse is not true: there are reduced sharbly
cycles that have size $>1$.  Nevertheless, experimentally we found
that every cohomology class we considered can be represented by a
sharbly cycle of size $1$.  Thus size gives a convenient measure of
success in our reduction algorithm.

\subsection{} We now turn to the reduction algorithm itself.  The
overall structure proceeds as described in \cite{experimental,
heckants}, and we refer to there for more details.  Let $\xi$ be a
nonreduced $1$-sharbly chain.  For each $1$-sharbly $\uu$ in the
support of $\xi $, we $\Gamma$-equivariantly choose a collection of
\emph{reducing points} for the nonreduced maximal faces of
$\uu$.\footnote{In \cite{experimental}, such points were called \emph{candidates}.}
These points, together with the original spanning points of the $\uu$,
are assembled into a new $1$-sharbly chain $\xi '$ homologous mod
$\Gamma$ to $\xi$ as described in \cite[Algorithm 4.13]{experimental}
and as recalled below.  Experimentally, the new chain $\xi '$ always
satisfies $\Size (\xi ') < \Size (\xi)$, although we cannot prove
this.  The reduction process is repeated until we obtain a $1$-sharbly
chain of size 1.

We now describe how individual $1$-sharblies $\uu$ are modified during
the algorithm.  Each $\uu$ is treated according to the number of its
nonreduced faces (cf.~Figure~\ref{fig:red}).  Suppose $\uu = [x_1, x_2, x_3,
x_4]$, and for any $i=1,\dotsc ,4$ let $f_i$ denote the $0$-sharbly
face $[\dots, \hat{x}_i, \dots]$.

Suppose first that  $\uu$ has four
nonreduced faces.  Then for each 
face $f_i$, we choose a reducing point $w_i$.  The reduction step
replaces $\uu$ with a sum of eleven new $1$-sharblies
\begin{equation} \label{eq:red4}
\begin{split} 
\uu \longmapsto &
-[x_1, x_4, w_2, w_3] +  [x_1, w_2, w_3, w_4] + [x_1, x_3, w_2, w_4] +
[x_2,x_4, w_1,w_3] \\
&-[x_2, w_1 ,w_3, w_4] - [x_2, x_3, w_1, w_4] - [x_1, x_2, w_3, w_4] -
[x_3, x_4, w_1, w_2]\\&-[x_4, w_1, w_2, w_3] - [w_1, w_2, w_3, w_4] +
[x_3, w_1, w_2, w_4].
\end{split}
\end{equation}

If $\uu$ has three nonreduced faces, arrange that the face $f_4$ is
reduced, and choose reducing points $w_1$, $w_2$, and $w_3$. Then the
reduction step replaces $\uu$ with a sum of eight new $1$-sharblies
\begin{equation} \label{eq:red3}
\begin{split}
\uu \longmapsto &
[x_1, x_2, x_3, w_1] - [x_1, x_3, w_1, w_2] + [x_1, w_1, w_2, w_3]
 + [x_1, x_2, w_1, w_3] \\
&- [x_3, x_4, w_1, w_2] - [x_4, w_1, w_2, w_3] + [x_2, x_4, w_1, w_3]
 - [x_1, x_4, w_2, w_3]. 
\end{split}
\end{equation}

If $\uu$ has two nonreduced faces, arrange that $f_1$ and $f_2$ are
not reduced, and choose reducing points $w_1$ and $w_2$.  Then the
reduction step replaces $\uu$ with a sum of five new 
$1$-sharblies 
\begin{equation} \label{eq:red2}
\begin{split}
\uu \longmapsto &
-[x_1, x_2, x_4, w_1] - [x_1, x_3, w_1, w_2] + [x_1, x_2, x_3, w_1]  +
[x_1, x_4, w_1, w_2] \\ &- [x_3, x_4, w_1, w_2].
\end{split}
\end{equation}

Finally, if $\uu$ has one nonreduced face, arrange that $f_1$ is not 
reduced, and choose reducing point $w_1$.  Then the reduction step
replaces $\uu$ with a sum of three new $1$-sharblies 
\begin{equation} \label{eq:red1}
\begin{split}
\uu \longmapsto &
[x_1, x_2, x_3, w_1] - [x_1, x_2, x_4, w_1] +  [x_1, x_3, x_4, w_1].
\end{split}
\end{equation}
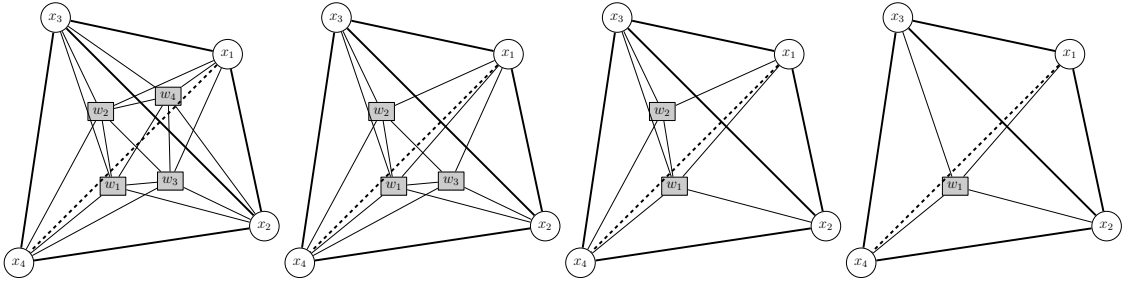
\begin{figure}
\resizebox{\textwidth}{!}{%
\tikzstyle{dot}=[circle,draw,minimum size=1mm]
\tikzstyle{redot}=[rectangle, draw, fill=black!20, minimum size=1mm]
\pgfmathsetmacro{\tetralength}{3}
\pgfmathsetmacro{\innertetralength}{1}
\begin{tikzpicture}[line join=bevel, z = -5, thick]
\node[dot] (x1) at (\tetralength,\tetralength,\tetralength) {$x_1$};
\node[dot] (x2) at (\tetralength,-\tetralength,-\tetralength) {$x_2$};
\node[dot] (x3) at (-\tetralength,\tetralength,-\tetralength) {$x_3$};
\node[dot] (x4) at (-\tetralength,-\tetralength,\tetralength) {$x_4$};
\node[redot] (w1) at (-\innertetralength, -\innertetralength -0.5, -\innertetralength) {$w_1$};
\node[redot] (w2) at (-\innertetralength, \innertetralength, \innertetralength) {$w_2$};
\node[redot] (w3) at (\innertetralength, -\innertetralength, \innertetralength) {$w_3$};
\node[redot] (w4) at (\innertetralength - 0.5, \innertetralength, -\innertetralength-0.5) {$w_4$};
\draw (w1) -- (w2) -- (w3) -- (w4) -- (w1) -- (w3);
\draw (w2) -- (w4);
\draw (w1) -- (x2);
\draw (w1) -- (x3);
\draw (w1) -- (x4);
\draw (w2) -- (x4);
\draw (w2) -- (x3);
\draw (w2) -- (x1);
\draw (w3) -- (x1);
\draw (w3) -- (x2);
\draw (w3) -- (x4);
\draw (w4) -- (x1);
\draw (w4) -- (x2);
\draw (w4) -- (x3);
\begin{scope}[ultra thick]
\draw (x4) -- (x2) -- (x3) -- (x4);
\draw (x1) -- (x2);
\draw (x1) -- (x3); 
\draw [dashed] (x1) -- (x4);
\end{scope}
\end{tikzpicture}

% reduce 3
\begin{tikzpicture}[line join=bevel, z = -5, thick]
\node[dot] (x1) at (\tetralength,\tetralength,\tetralength) {$x_1$};
\node[dot] (x2) at (\tetralength,-\tetralength,-\tetralength) {$x_2$};
\node[dot] (x3) at (-\tetralength,\tetralength,-\tetralength) {$x_3$};
\node[dot] (x4) at (-\tetralength,-\tetralength,\tetralength) {$x_4$};
\node[redot] (w1) at (-\innertetralength, -\innertetralength -0.5, -\innertetralength) {$w_1$};
\node[redot] (w2) at (-\innertetralength, \innertetralength, \innertetralength) {$w_2$};
\node[redot] (w3) at (\innertetralength, -\innertetralength, \innertetralength) {$w_3$};
\draw (x1) -- (w1);
\draw (x2) -- (w1);
\draw (x3) -- (w1);
\draw (w1) -- (w2) -- (w3) -- (w1);
\draw (x1) -- (w2);
\draw (x3) -- (w2);
\draw (w3) -- (x1);
\draw (w3) -- (x2);
\draw (w1) -- (x4);
\draw (w2) -- (x4);
\draw (w3) -- (x4);
\begin{scope}[ultra thick]
\draw (x4) -- (x2) -- (x3) -- (x4);
\draw (x1) -- (x2);
\draw (x1) -- (x3); 
\draw [dashed] (x1) -- (x4);
\end{scope}
\end{tikzpicture}

% reduce 2
\begin{tikzpicture}[line join=bevel, z = -5, thick]
\node[dot] (x1) at (\tetralength,\tetralength,\tetralength) {$x_1$};
\node[dot] (x2) at (\tetralength,-\tetralength,-\tetralength) {$x_2$};
\node[dot] (x3) at (-\tetralength,\tetralength,-\tetralength) {$x_3$};
\node[dot] (x4) at (-\tetralength,-\tetralength,\tetralength) {$x_4$};
\node[redot] (w1) at (-\innertetralength, -\innertetralength -0.5, -\innertetralength) {$w_1$};
\node[redot] (w2) at (-\innertetralength, \innertetralength, \innertetralength) {$w_2$};
\draw (x1) -- (w1);
\draw (x2) -- (w1);
\draw (x4) -- (w1);
\draw (w1) -- (w2);
\draw (w1) -- (x3);
\draw (w2) -- (x3);
\draw (w2) -- (x1);
\draw (w2) -- (x4);
\begin{scope}[ultra thick]
\draw (x4) -- (x2) -- (x3) -- (x4);
\draw (x1) -- (x2);
\draw (x1) -- (x3); 
\draw [dashed] (x1) -- (x4);
\end{scope}
\end{tikzpicture}

% reduce 1
\begin{tikzpicture}[line join=bevel, z = -5, thick]
\node[dot] (x1) at (\tetralength,\tetralength,\tetralength) {$x_1$};
\node[dot] (x2) at (\tetralength,-\tetralength,-\tetralength) {$x_2$};
\node[dot] (x3) at (-\tetralength,\tetralength,-\tetralength) {$x_3$};
\node[dot] (x4) at (-\tetralength,-\tetralength,\tetralength) {$x_4$};
\node[redot] (w1) at (-\innertetralength, -\innertetralength -0.5, -\innertetralength) {$w_1$};
\draw (x1) -- (w1);
\draw (x2) -- (w1);
\draw (x3) -- (w1);
\draw (x4) -- (w1);
\begin{scope}[ultra thick]
\draw (x4) -- (x2) -- (x3) -- (x4);
\draw (x1) -- (x2);
\draw (x1) -- (x3); 
\draw [dashed] (x1) -- (x4);
\end{scope}
\end{tikzpicture}
}
\caption{From left to right, we have from $4$ to $1$ nonreduced faces,
 as described
  in \eqref{eq:red4}, \eqref{eq:red3}, \eqref{eq:red2}, and
  \eqref{eq:red1}.} \label{fig:red} 
\end{figure}

\section{Computational results }\label{s:results}

\subsection{Overview of the data}\label{subsec:overview} We computed
$H^{5} (\Gamma_{0} (\n); \C)$ for ideals $\n \subset \OF$ of norm less
than or equal to $919$.  If a level $\n$ satisfies $\n \not = \bar
\n$, then we only computed the cohomology for one of $\n , \bar \n$.
There are $515$ ideals of norm less than or equal to $919$ which yield
$303$ up to conjugation.  Of these, $108$ levels had nontrivial
cohomology.  For each of these levels, we give relevant data in
Table~\ref{tab:nontrivial}.  We give the dimension $d_3$ of $H^{5}
(\Gamma_{0} (\n); \C)$, the dimension of the new subspace $d_3^{\new}$
as described in \S\ref{subsec:old}.  We give the dimension
$c_2^{\new}$ of the new cuspidal subspace of Bianchi forms as
described in \S\ref{subsec:gl2}.  We also give the number $g^{\prim}$
of primitive Hecke Gr\"ossencharacters as described in
\S\ref{subsec:grossen}. The last column $\Delta$ is the dimension of
the space that remains.

The ideals are indexed using an HNF-label, which encodes their Hermite
normal forms with respect to a certain basis of $\OF$.  For any ideal
$\n \subset \OF$, there are rational integers $a, c, d$ such that
$\set{a,d\w + c}$ is a $\ZZ$-basis of $\n$, with $0 \leq c < a$;
$d>0$; $ad = \Norm(\n)$; $d \mid a$; and $d \mid c$. The HNF-label is
the uniquely determined triple $[N, c, d]$.

\begin{center}
\begin{longtable}{llrrrrr}
\caption{Data for nontrivial cohomology groups $H^{5} (\Gamma_{0}
  (\n); \C)$ for $\Norm(\n) \leq 919$.  We give the dimension $d_3$ of $H^{5}
(\Gamma_{0} (\n); \C)$, the 
dimension of the new subspace $d_3^{\new}$ as described in
\S\ref{subsec:old}, the 
dimension $c_2^{\new}$ of the new cuspidal subspace of 
Bianchi forms as described in \S\ref{subsec:gl2}, and the number
$g^{\prim}$ of primitive Hecke Gr\"ossencharacters as described 
in \S\ref{subsec:grossen}.  The last column $\Delta$ is the
dimension of the space that remains, and these forms are dicsussed in
\S\ref{subsec:delta}.}  
\label{tab:nontrivial} \\
\toprule
$\HNF(\n)$ & $\HNF(\n')$ & $d_3$ & $d_3^{\new}$ & $g^{\prim}$ & $c_2^{\new}$ & $\Delta$\\
\midrule
\endfirsthead

\multicolumn{7}{c}%
{{\bfseries \tablename\ \thetable{} -- continued from previous page}} \\
\toprule
$\HNF(\n)$ & $\HNF(\n')$ & $d_3$ & $d_3^{\new}$ & $g^{\prim}$ & $c_2^{\new}$ & $\Delta$\\
\midrule
\endhead

\midrule
\multicolumn{7}{c}{{Continued on next page}} \\ 
\bottomrule
\endfoot

\bottomrule
\endlastfoot
{[49,18,1]} & {[49,30,1]} & 2 & 2 & 1 & 0 & 0\\
{[73,8,1]} & {[73,64,1]} & 2 & 2 & 0 & 1 & 0\\
{[75,5,5]} & {[75,5,5]} & 2 & 2 & 0 & 1 & 0\\
{[81,0,9]} & {[81,0,9]} & 2 & 2 & 1 & 0 & 0\\
{[121,0,11]} & {[121,0,11]} & 2 & 2 & 0 & 1 & 0\\
{[124,10,2]} & {[124,50,2]} & 2 & 2 & 0 & 1 & 0\\
{[144,0,12]} & {[144,0,12]} & 2 & 2 & 1 & 0 & 0\\
{[147,7,7]} & {[147,7,7]} & 2 & 2 & 0 & 1 & 0\\
{[147,67,1]} & {[147,79,1]} & 6 & 0 & 0 & 0 & 0\\
{[169,22,1]} & {[169,146,1]} & 4 & 4 & 2 & 0 & 0\\
{[171,21,3]} & {[171,33,3]} & 2 & 2 & 0 & 1 & 0\\
{[192,8,8]} & {[192,8,8]} & 2 & 2 & 0 & 1 & 0\\
{[196,0,14]} & {[196,0,14]} & 2 & 2 & 0 & 1 & 0\\
{[196,36,2]} & {[196,60,2]} & 6 & 0 & 0 & 0 & 0\\
{[219,64,1]} & {[219,154,1]} & 6 & 0 & 0 & 0 & 0\\
{[225,0,15]} & {[225,0,15]} & 6 & 0 & 0 & 0 & 0\\
{[228,14,2]} & {[228,98,2]} & 2 & 2 & 0 & 1 & 0\\
{[241,15,1]} & {[241,225,1]} & 2 & 2 & 0 & 1 & 0\\
{[243,9,9]} & {[243,9,9]} & 6 & 0 & 0 & 0 & 0\\
{[256,0,16]} & {[256,0,16]} & 4 & 4 & 2 & 0 & 0\\
{[273,16,1]} & {[273,256,1]} & 2 & 2 & 0 & 1 & 0\\
{[273,100,1]} & {[273,172,1]} & 2 & 2 & 0 & 1 & 0\\
{[283,44,1]} & {[283,238,1]} & 2 & 2 & 0 & 1 & 0\\
{[289,0,17]} & {[289,0,17]} & 2 & 2 & 0 & 1 & 0\\
{[292,16,2]} & {[292,128,2]} & 6 & 0 & 0 & 0 & 0\\
{[300,10,10]} & {[300,10,10]} & 8 & 2 & 0 & 1 & 0\\
{[324,0,18]} & {[324,0,18]} & 8 & 2 & 0 & 1 & 0\\
{[343,14,7]} & {[343,28,7]} & 8 & 2 & 0 & 1 & 0\\
{[343,18,1]} & {[343,324,1]} & 6 & 0 & 0 & 0 & 0\\
{[361,0,19]} & {[361,0,19]} & 2 & 2 & 0 & 1 & 0\\
{[361,68,1]} & {[361,292,1]} & 6 & 6 & 3 & 0 & 0\\
{[363,11,11]} & {[363,11,11]} & 8 & 2 & 0 & 1 & 0\\
{[372,50,2]} & {[372,134,2]} & 6 & 0 & 0 & 0 & 0\\
{[379,51,1]} & {[379,327,1]} & 2 & 2 & 0 & 1 & 0\\
{[399,163,1]} & {[399,235,1]} & 2 & 2 & 0 & 1 & 0\\
{[400,0,20]} & {[400,0,20]} & 2 & 2 & 0 & 1 & 0\\
{[412,92,2]} & {[412,112,2]} & 2 & 2 & 0 & 1 & 0\\
{[417,181,1]} & {[417,235,1]} & 2 & 2 & 0 & 1 & 0\\
{[432,12,12]} & {[432,12,12]} & 6 & 0 & 0 & 0 & 0\\
{[441,0,21]} & {[441,0,21]} & 8 & 2 & 0 & 1 & 0\\
{[441,54,3]} & {[441,90,3]} & 14 & 2 & 1 & 0 & 0\\
{[475,35,5]} & {[475,55,5]} & 2 & 2 & 0 & 1 & 0\\
{[481,211,1]} & {[481,269,1]} & 2 & 2 & 0 & 1 & 0\\
{[484,0,22]} & {[484,0,22]} & 6 & 0 & 0 & 0 & 0\\
{[496,20,4]} & {[496,100,4]} & 6 & 0 & 0 & 0 & 0\\
{[507,13,13]} & {[507,13,13]} & 6 & 6 & 0 & 3 & 0\\
{[507,22,1]} & {[507,484,1]} & 12 & 0 & 0 & 0 & 0\\
{[511,81,1]} & {[511,429,1]} & 6 & 0 & 0 & 0 & 0\\
{[511,137,1]} & {[511,373,1]} & 6 & 0 & 0 & 0 & 0\\
{[513,21,3]} & {[513,147,3]} & 8 & 2 & 0 & 1 & 0\\
{[523,60,1]} & {[523,462,1]} & 2 & 2 & 0 & 1 & 0\\
{[525,20,5]} & {[525,80,5]} & 6 & 0 & 0 & 0 & 0\\
{[529,0,23]} & {[529,0,23]} & 4 & 4 & 0 & 2 & 0\\
{[532,60,2]} & {[532,204,2]} & 4 & 4 & 0 & 2 & 0\\
{[553,102,1]} & {[553,450,1]} & 2 & 2 & 0 & 1 & 0\\
{[567,18,9]} & {[567,36,9]} & 6 & 0 & 0 & 0 & 0\\
{[576,0,24]} & {[576,0,24]} & 12 & 0 & 0 & 0 & 0\\
{[577,213,1]} & {[577,363,1]} & 2 & 2 & 0 & 1 & 0\\
{[579,277,1]} & {[579,301,1]} & 4 & 4 & 0 & 2 & 0\\
{[588,14,14]} & {[588,14,14]} & 14 & 2 & 0 & 1 & 0\\
{[588,134,2]} & {[588,158,2]} & 18 & 0 & 0 & 0 & 0\\
{[603,87,3]} & {[603,111,3]} & 2 & 2 & 0 & 1 & 0\\
{[613,65,1]} & {[613,547,1]} & 2 & 2 & 0 & 1 & 0\\
{[625,0,25]} & {[625,0,25]} & 8 & 8 & 4 & 0 & 0\\
{[637,165,1]} & {[637,471,1]} & 6 & 0 & 0 & 0 & 0\\
{[637,263,1]} & {[637,373,1]} & 6 & 0 & 0 & 0 & 0\\
{[651,25,1]} & {[651,625,1]} & 2 & 2 & 0 & 1 & 0\\
{[651,67,1]} & {[651,583,1]} & 2 & 2 & 0 & 1 & 0\\
{[657,24,3]} & {[657,192,3]} & 12 & 0 & 0 & 0 & 0\\
{[673,255,1]} & {[673,417,1]} & 2 & 2 & 0 & 1 & 0\\
{[675,15,15]} & {[675,15,15]} & 12 & 0 & 0 & 0 & 0\\
{[676,0,26]} & {[676,0,26]} & 4 & 4 & 0 & 2 & 0\\
{[676,44,2]} & {[676,292,2]} & 12 & 0 & 0 & 0 & 0\\
{[679,158,1]} & {[679,520,1]} & 2 & 2 & 0 & 1 & 0\\
{[679,326,1]} & {[679,352,1]} & 2 & 2 & 0 & 1 & 0\\
{[684,42,6]} & {[684,66,6]} & 12 & 0 & 0 & 0 & 0\\
{[700,20,10]} & {[700,40,10]} & 2 & 2 & 0 & 1 & 0\\
{[703,26,1]} & {[703,676,1]} & 2 & 2 & 0 & 1 & 0\\
{[721,46,1]} & {[721,674,1]} & 2 & 2 & 0 & 1 & 0\\
{[723,256,1]} & {[723,466,1]} & 8 & 2 & 0 & 1 & 0\\
{[729,0,27]} & {[729,0,27]} & 19 & 7 & 2 & 1 & 1\\
{[739,320,1]} & {[739,418,1]} & 2 & 2 & 0 & 0 & 2\\
{[741,334,1]} & {[741,406,1]} & 2 & 2 & 0 & 1 & 0\\
{[756,24,6]} & {[756,96,6]} & 2 & 2 & 0 & 1 & 0\\
{[757,27,1]} & {[757,729,1]} & 2 & 2 & 0 & 1 & 0\\
{[768,16,16]} & {[768,16,16]} & 20 & 2 & 0 & 1 & 0\\
{[784,0,28]} & {[784,0,28]} & 6 & 0 & 0 & 0 & 0\\
{[784,72,4]} & {[784,120,4]} & 16 & 4 & 2 & 0 & 0\\
{[793,230,1]} & {[793,562,1]} & 2 & 2 & 0 & 1 & 0\\
{[804,74,2]} & {[804,326,2]} & 2 & 2 & 0 & 1 & 0\\
{[819,27,3]} & {[819,243,3]} & 6 & 0 & 0 & 0 & 0\\
{[819,48,3]} & {[819,222,3]} & 6 & 0 & 0 & 0 & 0\\
{[832,24,8]} & {[832,72,8]} & 4 & 4 & 0 & 2 & 0\\
{[837,75,3]} & {[837,201,3]} & 4 & 4 & 0 & 1 & 2\\
{[841,0,29]} & {[841,0,29]} & 4 & 4 & 0 & 2 & 0\\
{[847,22,11]} & {[847,44,11]} & 6 & 0 & 0 & 0 & 0\\
{[849,238,1]} & {[849,610,1]} & 6 & 0 & 0 & 0 & 0\\
{[853,220,1]} & {[853,632,1]} & 4 & 4 & 0 & 1 & 2\\
{[867,17,17]} & {[867,17,17]} & 14 & 8 & 0 & 3 & 2\\
{[868,50,2]} & {[868,382,2]} & 6 & 0 & 0 & 0 & 0\\
{[868,134,2]} & {[868,298,2]} & 12 & 6 & 0 & 3 & 0\\
{[871,230,1]} & {[871,640,1]} & 2 & 2 & 0 & 1 & 0\\
{[876,128,2]} & {[876,308,2]} & 20 & 2 & 0 & 1 & 0\\
{[900,0,30]} & {[900,0,30]} & 26 & 2 & 0 & 1 & 0\\
{[903,79,1]} & {[903,823,1]} & 2 & 2 & 0 & 1 & 0\\
{[903,436,1]} & {[903,466,1]} & 2 & 2 & 0 & 1 & 0\\
{[912,28,4]} & {[912,196,4]} & 6 & 0 & 0 & 0 & 0\\
{[919,52,1]} & {[919,866,1]} & 4 & 4 & 0 & 2 & 0\\
\end{longtable}
\end{center}

\subsection{Eisenstein cohomology and interior cohomology}\label{subsec:eis}
Before presenting interpretations of our data, we recall some
background about Eisenstein cohomology \cite{harder}.
Recall that $X = \GL_3 (\C)/ \U (3)$ is the global symmetric space.  Let
$X^{\BorelSerre}$ be the partial compactification constructed by Borel
and Serre \cite{borel.serre}.  The quotient $Y := \Gamma _{0} (\n )\backslash X$
is an orbifold, and the quotient $Y^{\BorelSerre } := \Gamma_{0}
(\n )\backslash X^{\BorelSerre}$ is a compact orbifold with corners.  We
have 
\[
H^{*} (\Gamma_{0} (\n; \C) \simeq H^{*} (Y; \C) \simeq H^{*}
(Y^{\BorelSerre}; \C).
\]

Let $\partial Y^{\BorelSerre} = Y^{\BorelSerre}\smallsetminus Y$.  The
Hecke operators act on the cohomology of the boundary $H^{*} (\partial
Y^{\BorelSerre}; \C)$, and the inclusion of the boundary $\iota \colon
\partial Y^{\BorelSerre} \rightarrow Y^{\BorelSerre}$ induces a map on
cohomology $\iota^{*}\colon H^{*} (Y^{\BorelSerre}; \C) \rightarrow
H^{*} (\partial Y^{\BorelSerre}; \C)$ compatible with the Hecke
action.  The kernel $H^{*}_{!}  (Y^{\BorelSerre}; \C)$ of $\iota^{*}$
is called the \emph{interior cohomology}; it equals the image of the
cohomology with compact supports.  The goal of Eisenstein cohomology
is to use Eisenstein series and cohomology classes on the boundary to
construct a Hecke-equivariant section $s\colon H^{*} (\partial
Y^{\BorelSerre}; \C) \rightarrow H^{*} (Y^{\BorelSerre}; \C)$ mapping
onto a complement $H^{*}_{\Eis} (Y^{\BorelSerre}; \C )$ of the
interior cohomology in the full cohomology.  We call classes in the
image of $s$ \emph{Eisenstein classes}.

\subsection{Eisenstein cohomology I: classes coming from cuspforms on
\texorpdfstring{$\GL_2/F$}{GL2/F}} \label{subsec:gl2} Some cohomology
classes that we found correspond to weight $2$ Bianchi modular
cuspforms.  Such forms are the analogue of classical holomorphic
weight $2$ cuspforms for the group $\GL_{2} (\OF)$ \cite{egm}.  The
contribution of these forms is apparently via Eisenstein cohomology
\cite{harder} coming from the two types of maximal parabolic subgroups
of $\GL_{3} (F)$.  This phenomenon also occurs in the top degree
cohomology of subgroups of $\SL_{3} (\ZZ)$ \cite{agg}.

More precisely, let $f$ be a Bianchi cuspidal newform of level $\n$.
Then $f$ gives rise to $2$ eigenclasses $\phi$ and $\phi'$ in $H^{5}
(\Gamma_{0} (\n); \C)$ with associated Hecke polynomials
\begin{equation}\label{eq:gl2}
\begin{aligned}
\heckepoly(\phi, \p) &= (1-\Norm(\p)^2 t)(1  - a(\p,f) t + \Norm(\p) t^2)\\
\heckepoly(\phi', \p) &= (1 - t)(1 - \Norm(\p)a(\p,f)t + \Norm(\p)^3 t^2).  
\end{aligned}
\end{equation}
In particular, let $a(\p,i)$, $i = 1,2$ be the Hecke eigenvalues for
$\phi$, and let $a'(\p,i)$, $i = 1,2$ be the Hecke eigenvalues for
$\phi'$. Comparing to \eqref{eq:heckepoly}, we get 
\begin{equation}
\begin{aligned}
a(\p,1) &= a'(\p,2) = \Norm(\p)^2 + a(\p,f), \\
a(\p,2) &= a'(\p,1) =  1 + \Norm(\p) a(\p,f).
\end{aligned}
\end{equation}

\begin{example}\label{ex:73ec}
Let $\p_{73} = (-9\w + 1)$ be a prime ideal above $73$.  The space of
Bianchi cusp forms of level $\p_{73}$ is $1$-dimensional and new,
corresponding the isogeny class of the elliptic curve\footnote{This
and other similar labels refer to the \emph{$L$-functions and modular forms
database} \cite{lmfdb}.}
\curve{http://www.lmfdb.org/EllipticCurve/2.0.3.1/73.1/a/3}{73.1-a3}
with Weierstrass equation
\[y^2 + (\w + 1) x y + y = x^{3} + x^{2}.\]
This class contributes to the cohomology for $\GL_3$ over $F$ in two
ways with Hecke polynomials as given in \eqref{eq:gl2}.  The cohomology $H^{5}
(\Gamma_{0} (\p_{73}); \C)$ is
$2$-dimensional, so the Bianchi cuspform accounts for all of the 
cohomology at level $\p_{73}$.
\end{example}

\subsection{Eisenstein cohomology II: Hecke Gr\"ossencharacters}
\label{subsec:grossen} Some of the cohomology that arises can be
explained by Hecke Gr\"ossencharacters of weight $1$.  We expect that
these classes are also Eisenstein cohomology classes, again attached
to the two types of maximal parabolic subgroups in $\GL_{3}/F$.  We
briefly recall some of the details about Gr\"ossencharacters from
\cite{watkins}, specialized to our field $F$.

An $\infty$-type $\T$ is a pair of integers $[n_1, n_2]$ indexed by
the embeddings of $F$ in $\CC$.  The evaluation of $\T$ on an element
$\alpha \in F$ is given by $\T(\alpha) = \alpha_1^{n_1}
\alpha_2^{n_2}$.  For example, the $\infty$-type $[1,1]$ gives the
norm.  The \emph{weight} of $\T$ is the sum $n_1 + n_2$.  For a given
modulus $\m$, the $\infty$-type $\T$ is \emph{coherent} if $\T(\alpha)
= 1$ for all units $\alpha \equiv 1 \bmod{\m}$.  Given a modulus $\m$
and a coherent $\infty$-type $\T$, define $\G((\alpha)) = \T(\alpha)$
for all $\alpha \equiv 1 \bmod{\n}$.  This defines a Hecke
Gr\"ossencharacter $\G$ up to a finite quotient that is exactly the
Hecke character group for $\m$, the dual of the ray class group for
modulus $\m$.  \paulcomment{what is this?} In
particular, a Hecke character $\chi$ and an $\infty$-type $\T$
determine a Hecke Gr\"ossencharacter, which we denote $\G_{\chi,\T}$.

We observe primitive weight $1$ Hecke Gr\"ossencharacters contributing to the
cohomology as follows.  Let $H$ denote the Hecke character 
group of modulus $\m$, and let $D$ denote the Dirichlet character
group of modulus $\m$.  Suppose $\# H = \varphi(\#D)/2$ and $\#D \equiv 0
\pmod{6}$.  Let $\zeta$ be a primitive $\#D$th root of unity.  The values of
the Gr\"ossencharacters are in a subfield $K$ of $\QQ(\zeta)$ and give
rise to Hecke polynomials with coefficients in a real subfield $L$ of
$K$.  Note that $K$ is the cyclotomic field that is the codomain for
the Dirichlet characters, so $K$ is smaller than $\QQ(\zeta)$ in some
cases.  Let $\chi \in H$ be a character with primitive Gr\"ossencharacter
$\G_{\chi,[1,0]}$. Then $\chi$ contributes to the cohomology in two 
ways.  Specifically, $\chi$ gives rise to $2$ eigenclasses $\phi,
\phi'$ of level $\n$ with corresponding Hecke polynomials that have
coefficients in $L$ and factor over $K$ as 
\begin{equation}\label{eq:gross-factor}
\begin{aligned} 
\heckepoly(\phi,\p) &= (1 - \Norm(\p)^2 t)(1 - \G_{\chi,[1,0]}(\p) t)(1 -
\G_{\chi^{-1},[0,1]}(\p) t),\\
\heckepoly(\phi',\p) &= (1 - t)(1 - \Norm(\p)\G_{\chi,[1,0]}(\p)t)(1 -
\Norm(\p)\G_{\chi^{-1},[0,1]}(\p)t).
\end{aligned}
\end{equation}
Note $\G_{\chi^{-1},[0,1]}(\p)$ is the complex conjugate of
$\G_{\chi,[1,0]}(\p)$, so
 \[\G_{\chi,[1,0]}(\p)\G_{\chi^{-1},[0,1]}(\p) = \Norm(\p).\]
Let $a(\p,i)$, $i = 1,2$ be the Hecke eigenvalues for $\phi$, and let
$a'(\p,i)$, $i = 1, 2$ be the Hecke eigenvalues for $\phi'$.  
Comparing to \eqref{eq:heckepoly}, we find 
\begin{equation}
\begin{aligned}
a(\p,1) &= a'(\p,2)  = \Norm(\p)^2 + \G_{\chi,[1,0]}(\p) +
\G_{\chi^{-1},[0,1]}(\p)\\
a(\p,2) &= a'(\p,1) = 1 + \Norm(\p)(\G_{\chi,[1,0]}(\p) +
\G_{\chi^{-1},[0,1]}(\p)).
\end{aligned}  
\end{equation}
As in \S\ref{subsec:gl2}, we expect the classes
\eqref{eq:gross-factor} to be Eisenstein cohomology classes attached
to the maximal parabolic subgroups. 

We note that sometimes these classes correspond to elliptic curves $E$
over $F$ with CM by an order in $F$.  
In this case $\G_{\chi,[1,0]}(\p) + \G_{\chi^{-1},[0,1]}(\p)$ is in $\QQ$ for
all prime ideals $\p$.   The elliptic curve $E$ is a twist of a curve
$E'$ defined over $\Q$.  The base-change from $\Q$
to $F$ of the classical holomorphic cuspform corresponding to $E'$  is
not a cuspidal Bianchi modular form.
In fact, we have $L(E,s) = L(f,s)$, where $f$ is an Eisenstein
series.  Note in particular that such Bianchi modular forms apparently
contribute to the Eisenstein cohomology as in \S \ref{subsec:gl2}, although they are not cuspidal.

\begin{example} \label{ex:49CM}
Let $\p_7$ be a prime ideal above $7$, and let $\n =
\p_7^2$.  The cohomology $H^5(\Gamma_0(\n);\C)$ is $2$-dimensional
and can be explained by a Hecke Gr\"ossencharacter giving rise to a
elliptic curve over $F$ with CM by an order in $F$.

 The Hecke character group of modulus $\p_7$ is trivial.\footnote{We
used M.~Watkins's implementation of Hecke Gr\"ossencharacters in
\magma \cite{magma} in these computations.}  Let $e$ denote the
trivial character.  Then $e$ is not primitive, but the corresponding
Gr\"ossencharacter $\G_{e,[1,0]}$ is primitive.  The Dirichlet
character group of modulus $\p_7$ is cyclic of order $6$, so we are in
the situation described above with $\zeta$ a 6th root of unity.  Thus
$K = \field{http://www.lmfdb.org/NumberField/2.0.3.1}{\Q(\zeta)} = F$
and $L = \Q$.  In particular, the Gr\"ossencharacter gives rise to a
pair of Hecke polynomials that have rational coefficients and split
over $F$ as given in \eqref{eq:gross-factor}.  Because the Hecke
polynomials have $\Q$ coefficients, there should be an elliptic curve
over $F$ of conductor $\n$ with CM by an order in $F$.  A quick search
in LMFDB \cite{lmfdb} yields the curve
\curve{http://www.lmfdb.org/EllipticCurve/2.0.3.1/49.3/CMa/1}{49.3-CMa1}
with Weierstrass equation
\[y^2 + (\w + 1) y = x^{3} + (-\w - 1) x^{2} + \w
x - \w.\]

Let $\p_7'$ denote the other prime ideal above $7$.  Then
the cohomology at 
level $\p_7'^2$ is also $2$-dimensional and corresponds to the
isogeny class of the CM elliptic curve
\curve{http://www.lmfdb.org/EllipticCurve/2.0.3.1/49.1/CMa/1}{49.1-CMa1}.  
with Weierstrass equation
\[y^2 + \w y = x^{3} + (\w + 1) x^{2} + \w x.\]
\end{example}

\begin{example}
We now give an example that does not correspond to an elliptic curve.
Let $\p_{13}$ be a prime ideal above $13$.  The Hecke character group
is cyclic of order $2$, which we denote $\set{e, \chi}$.  The
corresponding Gr\"ossencharacters $\G_{e,[1,0]}$ and $\G_{\chi,[1,0]}$
are both primitive.  The Dirichlet character group of modulus
$\p_{13}$ is cyclic of order $12$, so we are in the situation
described above with $\zeta$ a $12$th root of unity.  We have Hecke
polynomials with coefficients in the real quadratic field $L =
\field{http://www.lmfdb.org/NumberField/2.2.12.1}{\Q(\sqrt{3})}$ that
split over the degree $4$ cyclotomic field $K =
\field{http://www.lmfdb.org/NumberField/4.0.144.1}{\Q(\zeta)}$. The
cohomology at level $\n = \p_{13}^2$ is $4$-dimensional, so this
accounts for all the cohomology.  Specifically, the Hecke operators
$T(\p,1)$ and $T(\p,2)$ diagonalize over $\Q(\sqrt{3})$.  The
resulting Hecke polynomials have coefficients in $\Q(\sqrt{3})$ and
factor over $\Q(\zeta)$ as given in \eqref{eq:gross-factor}.
\end{example}

\subsection{Old cohomology classes} \label{subsec:old} We observe
``old'' phenomena in cohomology consistent with what happens for
$\GL_{3}/\QQ$, according to van~Geeman, et.~al. \cite{top}, who cite
Reeder \cite{reeder}.  In particular, suppose $f$ is a newform for
$\GL_3$ over $F$ at level $\m$.  Then we observe $f$ contributing to
the cohomology at level $\n$ with multiplicity $3$ if $\n/\m$ is
prime, with multiplicity $6$ if $\n/\m$ is a square of a prime, and
with multiplicity $9$ if $\n/\m$ is a product of distinct primes.  We
call the subspace complementary to the old space \emph{new}.

\begin{example}
Let $\p_3$ be the prime above $3$, and let $\p_7$ be a prime above
$7$.  Let $\n$ be the norm $147$ ideal $\n = \p_3 \p_7^2$.  The
cohomology at level $\n$ is $6$-dimensional.  The cohomology at level
$\p_7^2$ is new and $2$-dimensional (see Example~\ref{ex:49CM}).  It
contributes to $H^5(\Gamma_0(\n); \C)$ with multiplicity $3$ since
$\n/\p_7^2$ is prime, so all of the cohomology at level $\n$ is old.
\end{example}

\begin{example}
Let $\n$ be the norm $576$ ideal $\n = (24) = \p_3^2 \q_2^3$,
where $\p_3$ is the prime above $3$ and $\q_2$ is the prime above $2$.  
The cohomology at level $(12) = \p_3^2\q_2^2$ is $2$-dimensional and
new, corresponding to the isogeny class of the CM elliptic curve
\curve{http://www.lmfdb.org/EllipticCurve/2.0.3.1/144.1/CMa/1}{144.1-CMa1}
with Weierstrass equation $y^2 = x^{3} + 1$.
The cohomology at level $\p_3\q_2^3$ is $2$-dimensional and new, corresponding
to isogeny class of the elliptic curve 
\curve{http://www.lmfdb.org/EllipticCurve/2.0.3.1/192.1/a/1}{192.1-a1}
with Weierstrass equation 
\[y^2 = x^{3} + \w x^{2} + (11 \w - 6) x + 11 \w - 1.\]
Each of these spaces contribute with multiplicity $3$ to the
cohomology at level $\n$ since $\n/(\p_3^2\q_2^2)$ and
$\n/(\p_3\q_2^3)$ are both prime.  The cohomology at level $\n$ is
$12$-dimensional, so oldforms account for all of the cohomology at
level $\n$.
\end{example}

\begin{example}
Let $\n$ be the norm $588$ ideal $\n = \p_3 \q_2 \p_7^2$, where $\p_3$
is the prime above $3$, $\q_2$ is the prime above $2$, and $\p_7$ is a
prime above $7$. The cohomology at level $\p_7^2$ is new and
$2$-dimensional (see Example~\ref{ex:49CM}).  It contributes with
multiplicity $9$ to the cohomology at level $\n$ since $\n/\p_7^2$ is
the product of two distinct primes. The cohomology at level $\n$ is
$18$-dimensional, so oldforms account for all of the cohomology at
level $\n$.
\end{example}

\begin{example}
Let $\n$ be the norm $657$ ideal $\n = \p_3^2 \p_{73}$, where $\p_3$
is the prime above $3$ and $\p_{73}$ is a prime above $73$.  The
cohomology at level $\p_{73}$ is $2$-dimensional (see example
Example~\ref{ex:73ec}).  It contributes with multiplicity $6$ to the
cohomology at level $\n$, since $\n/\p_{73}$ is the square of a prime
ideal. The cohomology at level $\n$ is $12$-dimensional, so oldforms
account for all of the cohomology at level $\n$.
\end{example}

\subsection{} \label{subsec:delta}
The remaining cases concern levels where, after removing the classes
corresponding to Eisenstein cohomology and old classes as above, we still had
some eigenclasses left over.  We believe these classes are in the
interior cohomology.  There were two different phenomena we observed:

\begin{itemize}
\item For one level $\n = \p_3^6$, where $\p_3$ is the prime over $3$,
we found an excess 1-dimensional space.  We believe this class is the
symmetric square lift of a weight 2 Bianchi modular cuspform.  Some
details can be found in \S\ref{subsec:n729}.
\item For four remaining levels, we found cohomology classes that
appear to coming from nonselfdual automorphic representations.  We
give details in \S\ref{subsec:nsd} below.
\end{itemize}

\subsection{Interior cohomology I: A symmetric square
lift}\label{subsec:n729} Let $\n = \p_3^6$ be the ideal of norm $729$
generated by $27$.  The cohomology $H^5(\Gamma_0(\n);\C)$ is
$19$-dimensional.  This is the only cohomology group in the range of
computations with odd dimension.

First we describe the old classes in the cohomology.
The cohomology at level $\p_3^4$ is new and $2$-dimensional,
coming from the isogeny class of 
\curve{http://www.lmfdb.org/EllipticCurve/2.0.3.1/81.1/CMa/1}{81.1-CMa1}
with Weierstrass equation
\[y^2 + y = x^3,\]
as described in \S\ref{subsec:grossen}.  The cohomology at level
$\p_3^5$ is $6$-dimensional, so it is all old, accounted for by the
cohomology at level $\p_3^4$.  Thus there is a $12$-dimensional old
subspace in the cohomology at level $\n$, as described in 
\S\ref{subsec:old}.  

Next we consider Eisenstein classes.  The space of Bianchi new cuspforms
of level $\n$ is $1$-dimensional, and corresponds to the isogeny
class of the elliptic curve
\curve{http://www.lmfdb.org/EllipticCurve/2.0.3.1/729.1/a/1}{729.1-a1}
with Weierstrass equation
\begin{equation}\label{eq:ec729}
y^2 + x y + \w y = x^3 -  x^2 - 2 \w x + \w.
\end{equation}
This contributes a $2$-dimensional subspace to the $\GL_3$ over $F$
cohomology at level $\n$, as described in \S\ref{subsec:gl2}.  

The Hecke character group at level $\p_3^3$ is cyclic of order
$3$. Let $\chi \in H$ denote a generator.  The Dirichlet character
group is isomorphic to $C_3 \times C_6$, so we are in the situation
described in \S\ref{subsec:grossen} with $K = F$ and $L = \QQ$.  The
Gr\"ossencharacters $\bG_{\chi,[1,0]}$ and $\bG_{\chi^2,[1,0]}$ are
primitive.  They give rise to a $4$-dimensional subspace of the
cohomology at level $\n$ corresponding to the two isogeny classes of the
elliptic curves over $F$ of conductor $\n$ with CM in $\OF$. A search
in LMFDB \cite{lmfdb} gives 
\curve{http://www.lmfdb.org/EllipticCurve/2.0.3.1/729.1/CMa/1}{729.1-CMa1}
with Weierstrass equation
\[y^2 + \w y = x^3,\]
and
\curve{http://www.lmfdb.org/EllipticCurve/2.0.3.1/729.1/CMb/1}{729.1-CMb1}
with Weierstrass equation
\[y^2 + (\w + 1) y = x^3 - \w.\]

After accounting for the above eigenclasses, we find a $1$-dimensional
subspace left.  Let $\phi$ be the Hecke eigenclass spanning this
complement.  Some Hecke polynomials for $\phi$ are given in
Table~\ref{tab:n729}.

We claim that $\phi$ appears to be the symmetric square lift of the
class corresponding to a Bianchi newform
\curve{http://www.lmfdb.org/ModularForm/GL2/ImaginaryQuadratic/2.0.3.1/729.1/a/}{729.1-a}
of level $\n$, in fact of the same form attached to the elliptic curve
\eqref{eq:ec729} above.  Indeed, let $f$ be this form, and for any
prime $\p \nmid \n$ let $a_{\p}$ be the eigenvalue of the
$\GL_{2}$-Hecke operator $T_{\p}$.  Let $\alpha_{\p}$, $\beta_{\p}$ be
the complex roots of the Hecke polynomial $1-a_{\p} t + \Norm (\p)
t^{2}$.  Then one expects the Hecke polynomial of the symmetric square
$\Sym^{2} f$ to have the form
\begin{equation}\label{eq:symsq}
( 1 - \alpha_{\p}^{2}t) (1-\alpha_{\p} \beta_\p t) (1-\beta_\p^{2}t) =
1 - (a_{\p}^{2} - \Norm (\p )) T +  (a_{\p}^{2} - \Norm (\p )) \Norm
(\p) t^{2} - \Norm (\p)^{3}t^{3}.  
\end{equation}
The form $f$ has eigenvalues $a_{2} = -1$, $a_{\p_{7}} = a_{\p'_7} =
2$, and $a_{\p_{13}} = a_{\p '_{13}}= -1$, which after insertion into
\eqref{eq:symsq} recovers Table \ref{tab:n729}.\footnote{We remark
that there is more to be said about this Bianchi modular form $f$.  It
is itself a base change of a weight two newform $g$
(\curve{http://www.lmfdb.org/ModularForm/GL2/Q/holomorphic/81/2/1/a/}{81.2.1.a})
on $\Gamma_{0} (81) \subset \SL_{2} (\Z)$ with coefficients in the
quadratic field $\QQ (\eta)$ of discriminant $12$, where $\eta^{2}=3$.
The eigenvalue of $T_{p}$, $p\nmid 3$ on $g$ away is rational
(respectively in $\ZZ \cdot \eta$) exactly when $p$ is inert
(respectively splits) in our imaginary quadratic field $F$. Thus $g$
corresponds to an abelian surface with extra twist as in
\cite{cremona.fake}.}

\begin{table}
  \caption{Hecke polynomials for the symmetric square eigenclass at the norm $729$
    level $\p_3^6$, described in \S\ref{subsec:n729}.} \label{tab:n729}
\[\begin{array}{cr}
\toprule
 \HNF(\p) & \heckepoly(\phi,\p)
\\ \midrule
{[4,0,2]}
& (-4t+1)(16t^2+7t+1)
\\ \midrule
{[7,2,1]}
& (-7t+1)(49t^2+10t+1)
\\ \midrule
{[7,4,1]}
& (-7t+1)(49t^2+10t+1)
\\ \midrule
{[13,3,1]}
& (-13t+1)(169t^2+25t+1)
\\ \midrule
{[13,9,1]}
& (-13t+1)(169t^2+25t+1)
\\  \bottomrule
\end{array}
\]
\end{table}

\subsection{Interior cohomology II: Nonselfdual
classes}\label{subsec:nsd} We come at last to the most interesting
classes in the paper, those that apparently correspond to nonselfdual
automorphic representations.  For a recent reference on automorphic
forms and representations, we refer to \cite{g, gh1, gh2}.  For the
connections between automorphic forms and cohomology of arithmetic
groups, we refer to \cite{bw}.

Let $E$ be a number field, let $\A_{E}$ the adeles of $E$, and let
$\pi$ be a cuspidal automorphic representation for $\GL_{n} (\A_{E})$.
Then $\pi$ is called \emph{selfdual} if it is isomorphic to its
contragredient.  We say a cuspform $f$ is selfdual if its
corresponding automorphic representation is selfdual, and that a
cohomology class is selfdual if its corresponding automorphic form is.
We can similarly define \emph{nonselfdual} for these objects.

We return to $\GL_{3}$ over our field $F$.  Let $\p \nmid 3$ be a
prime ideal and suppose that $\phi$ is a Hecke eigenclass with Hecke polynomial
\begin{equation}\label{eq:heckepoly2}
\heckepoly(\phi,\p) = 1 - a(\p,1) t + a(\p,2)\Norm(\p)t^2 - \Norm(\p)^{3}t^3
\end{equation}
at $\p$.  If $\phi$ is not an Eisenstein class, then at least one
Hecke polynomial must be irreducible.  One also knows that the eigenvalues $a (\p,1)$ and $a (\p, 2)$
are complex conjugates:
\begin{equation}\label{eq:cconjs}
a (\p,1) = \overline{a (\p, 2)}.
\end{equation}
Moreover, if $\phi$ is in addition a nonselfdual cuspidal
class with Hecke polynomial \eqref{eq:heckepoly2}, then the field
determined by the Hecke eigenvalues should be nonreal, and there must be
another nonselfdual class $\phi'$ complementing  $\phi$ and with Hecke
polynomial at $\p$ given by 
\begin{equation}\label{eq:heckepoly3}
\heckepoly(\phi',\p) = 1 - a(\p,2) t + a(\p,1)\Norm(\p)t^2 - \Norm(\p)^{3}t^3.
\end{equation}

To summarize, when detecting nonselfdual cuspidal cohomology classes, one
looks for the following criteria:
\begin{itemize}
\item One wants to see (at least one) irreducible Hecke
polynomial.\footnote{Strictly speaking, this does not imply that the
cohomology class is cuspidal, only that it appears in the
\emph{interior} cohomology of $\Gamma$ (see \cite{harder, agg} for the
definition).  For $\GL_{3}/\QQ$ being interior implies cuspidal.  For
our purposes, we will abuse notation and ignore this distinction.}
\item The Hecke eigenvalues should determine a nonreal field.
\item The classes should come in pairs, with the pattern of
eigenvalues given by \eqref{eq:heckepoly2}--\eqref{eq:heckepoly3}. 
\end{itemize}

We observed classes with these properties at four levels having norms
$739$, $837$, $853$, and $867$.  The classes at level norm $867$
appear to be base changes of classes on $\Gamma_{0} (153)\subset
\SL_{3} (\ZZ)$ originally found in \cite{top}.  We give details in
Examples \ref{ex:n739}--\ref{ex:n867} below.

\begin{example} \label{ex:n739}
Let $\p_{739}$ be the prime ideal above $739$ generated by $7 \w - 30$.
The cohomology at level $\p_{739}$ is $2$-dimensional.  Since the level is
prime, there are no oldforms or forms coming from Hecke
Gr\"ossencharacters.  The space of Bianchi cuspforms is trivial at
level $\n$.  Thus the cohomology $H^5(\Gamma_0(\p_{739});\C)$ is different
from the types explained above.  

There are $2$ eigenclasses $\phi$ and $\phi'$ with corresponding Hecke
polynomials that have coefficients in the imaginary quadratic field
$\field{http://www.lmfdb.org/NumberField/2.0.123.1}{\Q(\sqrt{-123})}$. The
Hecke polynomials for $\phi$ and $\phi'$ are given in
Table~\ref{tab:n739}.

Let $a(\p,i)$, $i = 1,2$ be the Hecke eigenvalues for $\phi$, and let
$a'(\p,i)$, $i = 1, 2$ be the Hecke eigenvalues for $\phi'$.  Then we
observe that 
\begin{equation}
a(\p,1) = a'(\p,2) = \overline{a'(\p,1)} = \overline{a(\p,2)}
\end{equation}
in the range of the computation.

\begin{table}
  \caption{Hecke polynomials for the nonselfdual eigenclasses at level
    $\p_{739}$ described in Example~\ref{ex:n739}.  Here, $\alpha$
    generates $\Q(\sqrt{-123})$ and satisfies the polynomial $x^2 - x
    + 31$.} \label{tab:n739}   
\[\begin{array}{cr}
\toprule
\HNF(\p) & \text{$\heckepoly(\phi_1,\p)$ and $\heckepoly(\phi_2,\p)$}   
\\ \midrule
 {[3,1,1]} & 
-27t^3+(-3\alpha-9)t^2+(-\alpha+4)t+1 
\\ & -27t^3+(3\alpha-12)t^2+(\alpha+3)t+1
\\ \midrule
{[4,0,2]} & 
-64t^3+(4\alpha-12)t^2+(\alpha+2)t+1 
\\ & -64t^3+(-4\alpha-8)t^2+(-\alpha+3)t+1
\\ \midrule
{[7,2,1]} & 
-343t^3+(-7\alpha-14)t^2+(-\alpha+3)t+1 
\\ & -343t^3+(7\alpha-21)t^2+(\alpha+2)t+1
\\ \midrule
{[7,4,1]} & 
-343t^3+(7\alpha-28)t^2+(\alpha+3)t+1 
\\ & -343t^3+(-7\alpha-21)t^2+(-\alpha+4)t+1
\\ \midrule
{[13,3,1]} & 
-2197t^3+(13\alpha+78)t^2+(\alpha-7)t+1 
\\ & -2197t^3+(-13\alpha+91)t^2+(-\alpha-6)t+1
\\ \midrule
{[13,9,1]} & 
-2197t^3+(-13\alpha-78)t^2+(-\alpha+7)t+1 
\\ & -2197t^3+(13\alpha-91)t^2+(\alpha+6)t+1
\\ \midrule
{[19,7,1]} & 
-6859t^3+(57\alpha-171)t^2+(3\alpha+6)t+1 
\\ & -6859t^3+(-57\alpha-114)t^2+(-3\alpha+9)t+1
\\ \midrule
{[19,11,1]} & 
-6859t^3+(-57\alpha-152)t^2+(-3\alpha+11)t+1 
\\ & -6859t^3+(57\alpha-209)t^2+(3\alpha+8)t+1
\\ \midrule
{[25,0,5]} & 
-15625t^3+(-25\alpha-125)t^2+(-\alpha+6)t+1 
\\ & -15625t^3+(25\alpha-150)t^2+(\alpha+5)t+1
\\ \midrule
{[31,5,1]} & 
-29791t^3+(-31\alpha+1302)t^2+(-\alpha-41)t+1 
\\ & -29791t^3+(31\alpha+1271)t^2+(\alpha-42)t+1
\\ \midrule
{[31,25,1]} & 
-29791t^3+(-155\alpha-744)t^2+(-5\alpha+29)t+1 
\\ & -29791t^3+(155\alpha-899)t^2+(5\alpha+24)t+1
\\  \bottomrule
\end{array}
\]
\end{table}
\end{example}

\begin{example} \label{ex:n837}
Let $\n$ be the ideal of norm $837 = 3^3 \cdot 31$ generated by
$21 \w + 12$.  The cohomology at level $\n$ is $4$-dimensional.   

The space of Bianchi cuspforms is $1$-dimensional and new at level
$\n$, coming from the isogeny class of the elliptic curve
\curve{http://www.lmfdb.org/EllipticCurve/2.0.3.1/837.1/a/1}{837.1-a1}
with Weierstrass equation 
\[y^2 + x y + \w y = x^3 + \w x^2 -  x.\]
This
contributes a $2$-dimensional subspace to $H^5(\Gamma_0(\n);\C)$, as
described in \S\ref{subsec:gl2}.   

Let $\phi$ and $\phi'$ be Hecke eigenclasses that generate the remaining
$2$-dimensional subspace of the cohomology.  The corresponding Hecke
polynomials have coefficients in $F$ and are given in Table~\ref{tab:n837}. 

Let $a(\p,i)$, $i = 1,2$ be the Hecke eigenvalues for $\phi$, and let
$a'(\p,i)$, $i = 1, 2$ be the Hecke eigenvalues for $\phi'$.  Then we
observe that as before,
\begin{equation}
a(\p,1) = a'(\p,2) = \overline{a'(\p,1)} = \overline{a(\p,2)}
\end{equation}
in the range of the computation.

\begin{table} 
  \caption{Hecke polynomials for the nonselfdual eigenclasses at the norm
    $837$ level generated by $21\w + 12$  described in
    Example~\ref{ex:n837}.} \label{tab:n837} 
\[\begin{array}{cr}
\toprule
 \HNF(\p) & \text{$\heckepoly(\phi_1,\p)$ and $\heckepoly(\phi_2,\p)$}   
\\ \midrule
{[4,0,2]} &
-64t^3+(-24\w+4)t^2+(-6\w+5)t+1 
\\ & -64t^3+(24\w-20)t^2+(6\w-1)t+1
\\ \midrule
{[7,2,1]} &
(-7t+1)(49t^2+10t+1)
\\ & (-7t+1)(49t^2+10t+1)
\\ \midrule
{[7,4,1]} &
(7t+1)(-49t^2+(12\w-6)t+1)
\\ & (7t+1)(-49t^2+(-12\w+6)t+1)
\\ \midrule
{[13,9,1]} &
(-13t+1)(169t^2+4t+1)
\\ & (-13t+1)(169t^2+4t+1)
\\  \bottomrule
\end{array}
\]
\end{table}

\end{example}

\begin{example} \label{ex:n853}
Let $\p_{853}$ be the prime ideal above $853$ generated by $27\w - 31$.
The cohomology $H^5(\Gamma_0(\p_{853});\C)$ is $4$-dimensional.  Since
the level is 
prime, there are no oldforms or forms coming from Hecke
Gr\"ossencharacters.  The space of Bianchi cuspforms is
$1$-dimensional and new at
level $\p_{853}$, coming from the the isogeny class of elliptic curve 
\curve{http://www.lmfdb.org/EllipticCurve/2.0.3.1/853.1/a/1}{853.1-a1}
over $F$ with Weierstrass equation
\[y^2 + x y + (\w + 1) y = x^3 + (\w - 1) x^2
+ (-\w - 1) x - \w,\]
as described in \S\ref{subsec:gl2}.  

Let $\phi$ and $\phi'$ be Hecke eigenclasses that generate the remaining
$2$-dimensional subspace of the cohomology.  The corresponding Hecke
polynomials that have coefficients in the imaginary quadratic field 
\field{http://www.lmfdb.org/NumberField/2.0.31.1}{$\Q(\sqrt{-31})$}
and are given in Table~\ref{tab:n853}. 

Let $a(\p,i)$, $i = 1,2$ be the Hecke eigenvalues for $\phi$, and let
$a'(\p,i)$, $i = 1, 2$ be the Hecke eigenvalues for $\phi'$.  Then as
in Example~\ref{ex:n739}, we
observe that 
\begin{equation}
a(\p,1) = a'(\p,2) = \overline{a'(\p,1)} = \overline{a(\p,2)}
\end{equation}
in the range of the computation.

\begin{table}
  \caption{Hecke polynomials for the nonselfdual eigenclasses at level
    $\p_{853}$ described in Example~\ref{ex:n853}. Here, $\beta$
    generates $\Q(\sqrt{-31})$ and satisfies the polynomial $x^2 - x +
    8$.} \label{tab:n853}   
\[\begin{array}{cr}
\toprule
 \HNF(\p) & \text{$\heckepoly(\phi_1,\p)$ and $\heckepoly(\phi_2,\p)$}
\\ \midrule
{[3,1,1]} &
(-3t+1)(9t^2+5t+1) \\ &
(-3t+1)(9t^2+5t+1)
\\ \midrule
{[4,0,2]} &
(4t+1)(-16t^2+(-2\beta+1)t+1) \\ &
(4t+1)(-16t^2+(2\beta-1)t+1)
\\ \midrule
{[7,2,1]} &
-343t^3+(28\beta-56)t^2+(4\beta+4)t+1\\  &
-343t^3+(-28\beta-28)t^2+(-4\beta+8)t+1
\\ \midrule
{[7,4,1]} &
(-7t+1)(49t^2+10t+1)\\ &
(-7t+1)(49t^2+10t+1) 
\\ \midrule
{[13,3,1]} &
-2197t^3+(104\beta-195)t^2+(8\beta+7)t+1\\ &
-2197t^3+(-104\beta-91)t^2+(-8\beta+15)t+1 
\\ \midrule
{[13,9,1]} &
-2197t^3+(-52\beta+26)t^2+(-4\beta+2)t+1\\ &
-2197t^3+(52\beta-26)t^2+(4\beta-2)t+1 
\\ \midrule
{[19,7,1]} &
-6859t^3-152\beta t^2+(-8\beta+8)t+1\\ &
-6859t^3+(152\beta-152)t^2+8\beta t+1 
\\ \midrule
{[19,11,1]} &
(-19t+1)(361t^2+23t+1)\\ &
(-19t+1)(361t^2+23t+1) 
\\  \bottomrule
\end{array}
\]
\end{table}
\end{example}

\begin{example}\label{ex:n867}
Let $\n$ be the norm $867$ ideal $\n = \p_3 \q_{17}$, where $\p_3$ is
the prime above $13$ and $\q_{17}$ is the prime above $17$.  The
cohomology $H^5(\Gamma_0(\n);\C)$ is $14$-dimensional.  There is a
$6$-dimensional contribution of old forms, as described in
\S\ref{subsec:old}, coming from level the $2$-dimensional cohomology at
level $\q_{17}$, so the new subspace of $H^5(\Gamma_0(\n);\C)$ is
$8$-dimensional.  The space of Bianchi cuspforms at level $\n$ is
$5$-dimensional.  The new subspace is $3$-dimensional, including a
rational cuspidal newform corresponding to the isogeny class of 
\curve{http://www.lmfdb.org/EllipticCurve/2.0.3.1/867.1/a/1}{867.1-a1}
with Weierstrass equation
\[y^2 + y = x^3 + x^2 - 59 x - 196.\]  The new subspace of Bianchi
forms contributes a $6$-dimensional subspace in $H^5(\Gamma_0(\n);\C)$
as described in \S\ref{subsec:gl2}.  Let $\phi$ and $\phi'$ be the
Hecke eigenclasses that generate the remaining $2$-dimensional
subspace of the cohomology.  The corresponding Hecke polynomials have
coefficients in the imaginary quadratic field
$\field{http://www.lmfdb.org/NumberField/2.0.8.1}{\Q(\sqrt{-2})}$ and
are given in Table~\ref{tab:n867}.

Let $a(\p,i)$, $i = 1,2$ be the Hecke eigenvalues for $\phi$, and let
$a'(\p,i)$, $i = 1, 2$ be the Hecke eigenvalues for $\phi'$.  Then as
in Example~\ref{ex:n739}, we
observe that 
\begin{equation}
a(\p,1) = a'(\p,2) = \overline{a'(\p,1)} = \overline{a(\p,2)}
\end{equation}
in the range of the computation.  

Since our classes meet the criteria in \S \ref{subsec:nsd}, we expect
them to be nonselfdual cuspidal classes.  This appears
to be the case, and in fact they appear to be base changes
from $\QQ$ to $F$ of cohomology classes on $\Gamma_{0} (153) \subset
\SL_{3} (\ZZ)$ appearing in \cite{top}.  Indeed, if so one expects the
$L$-functions to satisfy
\[
L (g,s) = L (f, s) L (f \otimes \chi , s),
\]
where $g$ (respectively $f$) is the base change on $\GL_{3}/F$ (resp.,
the original automorphic form on $\GL_3/\QQ$) and $\chi$ is the
quadratic character corresponding to the extension $F/\QQ$.  In
particular, one expects the following identification of Hecke
eigenvalues away from the primes $3$ and $17$:
\begin{align}
a (\p , 1) &= a (p, 1) \quad \text{if $\p \mid  p$ and $p$ splits,}\label{eq:basechange1}\\
a (\p , 1) &= a (p,1)^{2} - 2 p a (p,2) \quad  \text{if $\p  = (p)$
and $p$ is inert.}\label{eq:basechange2}
\end{align}
Here we have written $a (\p , k)$ for the eigenvalues of the class on
$\GL_{3} (\OF)$ and $a (p, k)$ for the eigenvalues of the class on $\SL_{3} (\ZZ)$.
From \cite{top}, the relevant Hecke eigenvalues $a (p,1)$ are 
\[
a (2,1) = 1, \quad a (7, 1) = -6\gamma -3, \quad a (13,1) = -9\gamma -12, 
\]
which agrees with \eqref{eq:basechange1}--\eqref{eq:basechange2} since
$2$ is inert and $7, 13$ both split in $F$. 

We remark that the classes at level norms $739$, $837$, and $853$
\emph{cannot} be base changes from $\Q$, since they have different
eigenvalues at some conjugate pairs of primes of $\OF$.
  
\begin{table}
  \caption{Hecke polynomials for nonselfdual eigenclasses at level
$\p_3\q_{17}$ of norm $867$ described in Example~\ref{ex:n867}. Here,
$\gamma = \sqrt{-2}$.} \label{tab:n867}
\[\begin{array}{cr}
\toprule
 \HNF(\p) & \text{$\heckepoly(\phi_1,\p)$ and $\heckepoly(\phi_2,\p)$}
\\ \midrule
{[4,0,2]}&
-64t^{3}-12t^{2}+3t + 1\\&
-64t^{3}-12t^{2}+3t + 1\\
\midrule
{[7,2,1]}&
-343t^3+(42\gamma-21)t^2+(6\gamma+3)t+1\\&
-343t^3+(-42\gamma-21)t^2+(-6\gamma+3)t+1\\
\midrule
{[7,4,1]}&
-343t^3+(42\gamma-21)t^2+(6\gamma+3)t+1\\&
-343t^3+(-42\gamma-21)t^2+(-6\gamma+3)t+1\\
\midrule
{[13,9,1]}&
-2197t^3+(-156\gamma-117)t^2+(-12\gamma+9)t+1\\&
-2197t^3+(156\gamma-117)t^2+(12\gamma+9)t+1\\  
\bottomrule
\end{array}
\]
\end{table}
\end{example}
\bibliographystyle{amsplain_initials_eprint_doi_url}
\bibliography{gl3neg3}
\end{document}